\documentclass{article}

\usepackage[english]{babel}

\usepackage[letterpaper,top=2cm,bottom=2cm,left=3cm,right=3cm,marginparwidth=1.75cm]{geometry}

\usepackage{amsmath,amssymb,amsthm}
\usepackage{graphicx}
\usepackage[colorlinks=true, allcolors=blue]{hyperref}
\usepackage[capitalize]{cleveref} 
\usepackage[dvipsnames]{xcolor}

\theoremstyle{plain}
\newtheorem{axiom}{Axiom}
\newtheorem{claim}[axiom]{Claim}
\newtheorem{theorem}{Theorem}[section]
\newtheorem{lemma}[theorem]{Lemma}
\newtheorem{proposition}[theorem]{Proposition}
\newtheorem{conjecture}[theorem]{Conjecture}

\theoremstyle{definition}

\newtheorem{remark}[theorem]{Remark}

\usepackage{enumitem}  

\usepackage{todonotes}

\usepackage{enumitem} 

\newcommand{\bbZ}{\mathbb{Z}}
\newcommand{\N}{\mathbb{N}}
\newcommand{\E}{\mathbb{E}}
\renewcommand{\P}{\mathbb{P}}
\newcommand{\R}{\mathbb{R}}
\newcommand{\set}[1]{\left\{#1\right\}}

\newcommand{\eps}{\varepsilon}

\newcommand{\LIS}{L}
\newcommand{\LNS}{L^{(\leq)}}

\title{Longest increasing subsequences for distributions with atoms, and an inhomogeneous Hammersley process}
\author{Anne-Laure Basdevant\footnote{LPSM - UMR CNRS 8001, Sorbonne Universit\'e, 4 place Jussieu, 75005 Paris, France, {\it anne.laure.basdevant@normalesup.org}}, Lucas Gerin\footnote{CMAP, École Polytechnique, CNRS, Route de Saclay, 91128 Palaiseau Cedex, France
 {\it lucas.gerin@polytechnique.edu}}, Maxime Marivain\footnote{CMAP, École Polytechnique, CNRS, Route de Saclay, 91128 Palaiseau Cedex, France
 {\it maxime.marivain@polytechnique.edu}}}

\begin{document}
\maketitle

\begin{abstract}
We investigate the asymptotic behavior of the longest increasing subsequence $L_n$ in a  sequence of i.i.d. random variables, for an arbitrary distribution. If the common distribution is atomless then a famous result by Hammersley and Ver\v{s}ik-Kerov states that $L_n$ grows like $2\sqrt{n}$.

For any discrete distribution, we characterize the asymptotic order of $L_n$ through a variational problem and provide explicit estimates for classical distributions. The proofs rely on a coupling with an inhomogeneous version of the discrete-time continuous-space Hammersley process. This reveals that, in contrast to the continuous case, the discrete setting exhibits a wide range of growth rates between $\mathcal{O}(1)$ and $o(\sqrt{n})$, depending on the tail behavior of the distribution. We can then easily deduce the asymptotics of $L_n$ for a completely arbitrary distribution.
\end{abstract}

\medskip

\noindent{\bf Keywords:} Combinatorial probability, longest increasing subsequence,  Hammersley process.

\section{Introduction}
\subsection{Model and main results}
For $n \geq 1$ and real numbers $x_1, \dots, x_n$, denote by
$$
L(x_1,\dots,x_n)=\max\set{k\geq 1;\ \hbox{ there exists }1\leq i_1<\dots <i_k\leq n, \text{ such that }x_{i_1} < x_{i_2} < \dots < x_{i_k} }
$$
the length of one of the longest increasing subsequences in $x_1,\dots,x_n$.
The subject of interest of the present paper is the sequence $(\LIS_n)_{n\geq 1}$ of random variables defined by
$$
\LIS_n=L\left(X_1,\dots, X_n\right),
$$
where $(X_n)_{n\geq 1}$ is a sequence of i.i.d. random variables of common distribution $\rho$. 

If $\rho$ is atomless, then it is clear that $L(X_1,\dots, X_n)$
has the same distribution as the length $\mathcal{L}_n$ of the longest increasing subsequence in a uniformly random permutation of size $n$.
The study of $\mathcal{L}_n$ is often called the Hammersley problem, it is a historically important problem in combinatorial probability, not least because it can be tackled by several approaches combining algebraic combinatorics, calculus of variations, stochastic processes and particle systems, ... 
We briefly mention the main results concerning the asymptotics of $\mathcal{L}_n$, for more details and references we refer to \cite{Romik}. 
Hammersley \cite{Hamm} was the first to establish the existence of $c>0$ such that $\mathbb{E}[\mathcal{L}_n]=(c+o(1))\sqrt{n}$ and Ver\v{s}ik and Kerov \cite{vershik1977asymptotics} proved that $c=2$. 

In the present article, on the contrary, we deal with the case where the common distribution $\rho$ of $X_n$'s has atoms. 
We first present the results for the case of a distribution on $\mathbb{Z}_{\geq 1}$, and then show that the generic case presents no further difficulties.

We denote by $\mathbf{p}=(p_i)_{i\geq 1}$ the  probability mass function of $\rho$ when the support of $\rho$ is contained in $\mathbb{Z}_{\ge 1}$.
Let us first observe that if $U_1,\dots, U_n$ are i.i.d. uniform random variables in $(0,1)$, then 
$$
L(X_1,\dots, X_n)\leq L(X_1+U_1,\dots, X_n+U_n)\stackrel{\text{(d)}}{=} \mathcal{L}_n,
$$
so that $L_n= \mathcal{O}(\sqrt{n})$, in probability. We will see that the problem is quite rich in the discrete settings:
essentially all asymptotic behaviors between $\mathcal{O}(1)$  and $\mathrm{o}(\sqrt{n})$ are possible.
\begin{remark}
 If ${\bf p}$ has a finite support then it is very easy to check that 
$\LIS_n\to K$ a.s., where $K:=\mathrm{card}\left\{i,p_i>0\right\}$. Our main result  Theorem \ref{th:encadrement} is true in this case but immediate. Therefore we assume from now on that ${\bf p}$ does not have a bounded support and without loss of generality that $p_i>0$ for all $i\ge 1$.
\end{remark}

\begin{figure}
    \centering
    \includegraphics[width=13cm]{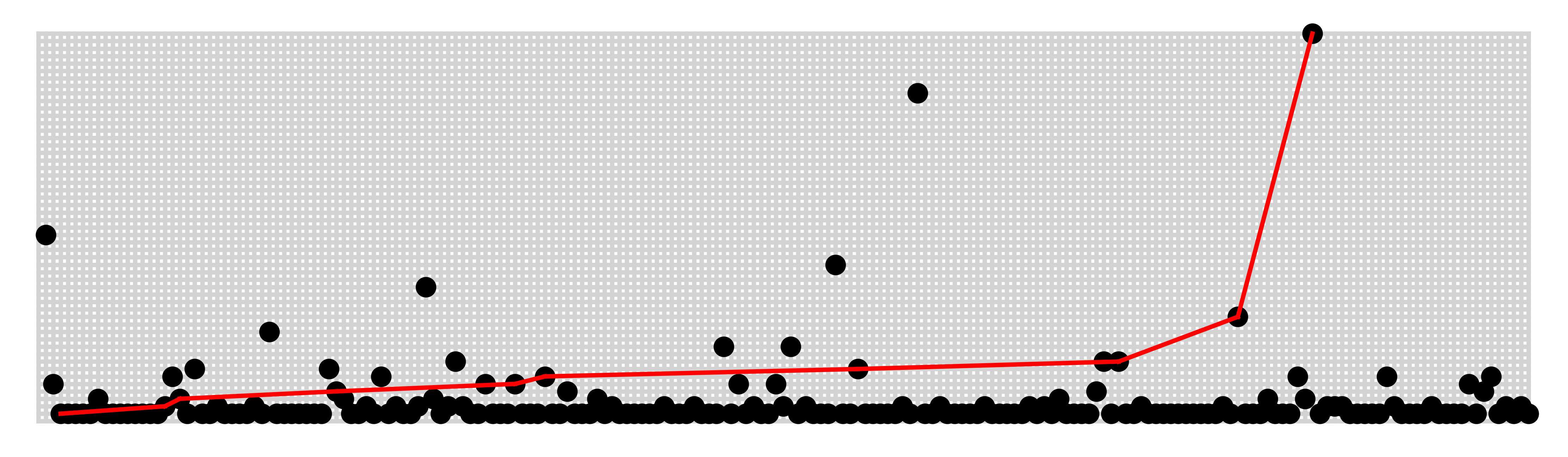}
    \caption{A sample of $X_1,\dots,X_n$ for $n=200$ and a power law $p_i\asymp i^{-2.2}$. In red: a longest increasing subsequence of size $\LIS_n=10$.}
    \label{fig:simu_puissance}
\end{figure}

\begin{figure}
    \centering
    \includegraphics[width=13cm]{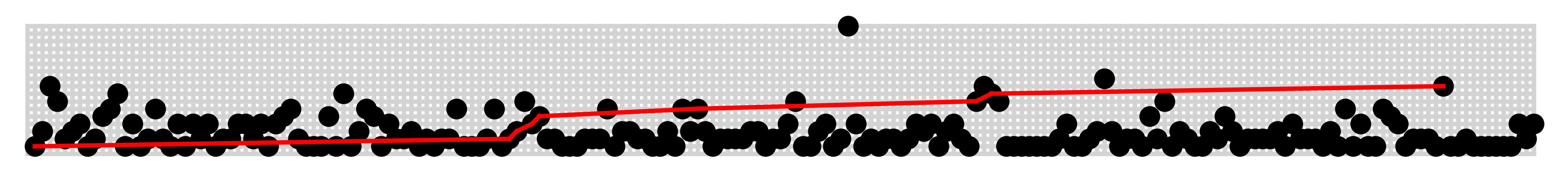}
    \caption{A sample of $X_1,\dots,X_n$ for $n=200$ and a geometric distribution  $p_i\asymp 0.6^i$. In red: a longest increasing subsequence of size $\LIS_n=9$.}
    \label{fig:simu_geo}
\end{figure}
We need a few notation in order to state our results.
For $t>0$, let $f_t$ be the solution of the variational problem
\begin{equation}\label{eq:fn}
    f_t= \inf_{x\geq 0}\left\{ x+\sum_{i\geq 1}\frac{tp_i}{tp_i+x}\right\}
\end{equation}
and 
let $w_t$ the unique positive solution of the equation  
\begin{equation}\label{eq:wn}
w_t=\sum_{i\geq 1} \frac{tp_i}{tp_i+w_t}.
\end{equation}
Here is our main theorem.
\begin{theorem}\label{th:encadrement}
Let $\mathbf{p}=(p_i)_{i\ge 1}$ be a discrete distribution on $\bbZ_{\geq 1}$, then
\begin{equation}\label{eq:Cv_ps_Theorem}
\lim_{n\to \infty}\frac{L_n}{\mathbb{E}[L_n]}= 1,\quad \text{a.s.}
\end{equation}
Furthermore
$$
\limsup_{n\to+\infty}\frac{L_n}{f_n}\leq 1\qquad \mbox{ and } \qquad 
\liminf_{n\to+\infty}\frac{L_n}{w_n}\geq 1 \qquad\mbox{a.s.}
$$
where $f_n,w_n$ are defined in  \cref{eq:fn,eq:wn}.
\end{theorem}

It is easy to see that $w_n\le f_n\le 2w_n$ (see Lemma \ref{lem:basic} \ref{encadref}). Thus Theorem \ref{th:encadrement} ensures that we have identified the right order of magnitude for $L_n$:
$$
\frac{1}{2}\leq \liminf_{n\to+\infty}\frac{L_n}{f_n}\leq \limsup_{n\to+\infty}\frac{L_n}{f_n}
\leq 1 \qquad \text{a.s.}
$$
We will also see in Lemma \ref{lem:basic} \ref{majorf} that $f_n=\mathrm{o}(\sqrt{n})$ for every ${\mathbf p}$, so that $L_n$ is always negligible with respect to the continuous case. 

At this stage, it is difficult to give an intuition behind formulas \cref{eq:fn,eq:wn}. These expressions will have a probabilistic interpretation in Section \ref{sec:Hammersley}, through couplings with a particle system. However, we can already mention that for several variants of Hammersley's problem, the constant factor in the first order of the asymptotics can be written as the solution of a variational problem. For instance a probabilistic proof of $c=2$ in the original Hammersley problem, due to Cator and Groeneboom \cite{CatorGroeneboom} is obtained by showing that $2=c=\inf_\lambda\{\lambda+\frac{1}{\lambda}\}$, which is somehow a continuous analogue of \cref{eq:fn} (see also \cite[Lem.2.2.]{seppalainen1997increasing} for a discrete Hammersley problem).

We will see below that for every $\mathbf{p}$ for which we obtain an equivalent for $\mathbb{E}[L_n]$, then $\mathbb{E}[L_n]\sim f_n$. This observation, together with the fact that the variational problem gives the right asymptotic constant for the above mentioned models, leads to the following conjecture.

\begin{conjecture}
For every  distribution $\mathbf{p}$ on $\bbZ_{\geq 1}$, it holds that  $\displaystyle{\frac{L_n}{f_n}} \stackrel{n\to+\infty }{\to}1$ a.s.
\end{conjecture}

Unfortunately it is not easy to find, for an arbitrary  $\mathbf{p}$, the asymptotic behavior of $w_n$ and $f_n$  directly from \cref{eq:fn} and \cref{eq:wn}. In the next results, we give the order of magnitude of $f_n$ in a variety of different cases.

Firstly, we give for comparison purposes the statement for Poisson and geometric distributions. It turns out that for these two cases we can directly determine the correct order of $\LIS_n$ without relying on Theorem \ref{th:encadrement}.

\begin{proposition}[Asymptotics of $L_n$: Poisson and geometric distributions]\label{prop:lighttail}
\ 
\begin{itemize}
\item[(i)] If $\mathbf{p}$ is a geometric distribution of parameter $p$, \emph{i.e.} $p_i=(1-p)^{i-1}p$ for $i\ge 1$ then
$$f_n\sim \frac{\log n}{|\log(1-p)|}.$$
\item[(ii)] 
If $\mathbf{p}$ is a Poisson distribution of parameter $\lambda$,  then
$$f_n\sim \frac{\log n}{\log(\log n)}. $$
\end{itemize}
Moreover, in both cases 
$$\lim_{n\to +\infty} \frac{\LIS_n}{f_n}=1 \mbox{ a.s.},
\qquad
\lim_{n\to +\infty} \frac{\mathbb{E}[\LIS_n]}{f_n}=1.
$$
\end{proposition}

A particularly interesting case in Theorem \ref{th:encadrement} is that of variables with heavy tails. In this case,  $f_n$ can be estimated and the asymptotic given in Theorem \ref{th:encadrement} does not seem to be obtainable by direct methods.

\begin{proposition}[Asymptotics of $L_n$: heavy-tailed distributions]\label{prop:heavytail} If the distribution function $\mathbf{p}$ is eventually non-increasing and:\  \begin{enumerate}[label=(\roman*)]
    \item if $p_i\stackrel{i\to+\infty}{\sim} c(\log i)^\gamma i^{-\beta}$ for some $\beta>1$ and $\gamma \in \R$, then, we have
$$ 0<\liminf_{n\to \infty} \frac{\LIS_n}{(n(\log n)^\gamma)^{\frac{1}{\beta+1}}}\le \limsup_{n\to \infty} \frac{\LIS_n}{(n(\log n)^\gamma)^{\frac{1}{\beta+1}}} <\infty
 \mbox{ a.s.}$$ \label{PropHT_i}
 \item if $p_i\stackrel{i\to+\infty}{\sim} c(\log i)^\gamma i^{-1}$ for some $\gamma<-1$, then, we have
$$ 0<\liminf_{n\to \infty} \frac{\LIS_n}{\sqrt{n(\log n)^{\gamma+1}}}\le \limsup_{n\to \infty} \frac{\LIS_n}{\sqrt{n(\log n)^{\gamma+1}}} <\infty
 \mbox{ a.s.}$$  \label{PropHT_ii}
\end{enumerate} 
The same results hold true if we replace    $L_n$ by $\E(L_n)$.
\end{proposition}

Proposition \ref{prop:heavytail} \ref{PropHT_i} implies that, as announced above, one can find for every $0<\xi<1/2$  and every $\zeta$ a distribution such that $L_n \asymp n^\xi(\log n)^\zeta$.
We give later in Lemma \ref{lem:Encadrement_f_n}, Proposition \ref{prop:f_sous_exp},  Lemma \ref{lem:borne_nu_mu} various  bounds for $f_t$ that apply to more general  distributions.

We can legitimately ask whether, for natural choices of ${\bf p}$, the asymptotic given by Theorem \ref{th:encadrement} can be obtained by more direct methods. Let us now briefly explain, in the case of heavy-tailed distributions, that the two bounds of Proposition \ref{prop:heavytail} differ from that obtained by two elementary strategies.
\begin{itemize}
    \item  Let us compare  Proposition \ref{prop:heavytail} with some naive bounds given by the following observation. Let $K_n=\mathrm{card}(\{X_1,X_2,\dots, X_n\})$ denote the number of distinct values in the sample, clearly $L_n\leq K_n$. On the other hand, the first occurrence of each distinct value induces a uniform random permutation of size $K_n$, so that with a little work one can show that for a general distribution ${\mathbf p}$, asymptotically with high probability,
$$
\sqrt{K_n}
\leq
\LIS_n(X_1,\dots, X_n)
\leq K_n.$$

When $p_i\stackrel{i\to+\infty}{\sim}  c(\log i)^\gamma i^{-\beta}$ with $\beta>1$, it is known that $K_n$ is of order $(\log n)^{\frac{\gamma}{\beta}}n^{1/\beta}$ (see for instance \cite[Cor.21]{Pitman}). Therefore, comparisons with $K_n$ would lead to bounds of the form
$$ \kappa\sqrt{(\log n)^{\frac{\gamma}{\beta}}n^{1/\beta}}\le L_n \le \kappa'(\log n)^{\frac{\gamma}{\beta}}n^{1/\beta}$$
 for some $\kappa,\kappa'>0$. Note that Proposition \ref{prop:heavytail} implies that in fact 
 $$ \sqrt{(\log n)^{\frac{\gamma}{\beta}}n^{1/\beta}}\ll L_n \ll c'(\log n)^{\frac{\gamma}{\beta}}n^{1/\beta}.$$
 Thus, in the case of heavy-tailed distributions, both the lower and upper bounds obtained thanks to Theorem \ref{th:encadrement} are sharper and nontrivial.
 \item We will describe in Section \ref{sec:greedy} a greedy algorithm whose output is an increasing subsequence. It provides a lower bound for $f_n$, which coincides with the lower bound of  Proposition \ref{prop:heavytail} \ref{PropHT_i}. However, in the case \ref{PropHT_ii} of Proposition \ref{prop:heavytail} (\emph{i.e.}, distributions of the form $p_i\stackrel{i\to+\infty}{\sim} c(\log i)^\gamma i^{-1}$)  the greedy algorithm fails to capture the correct asymptotic order. In this situation, the lower bound from Proposition \ref{prop:heavytail} becomes particularly relevant
 (see also the discussion in Remark \ref{rem:indice1}).
\end{itemize}
These two points seem to confirm that the particle system approach is necessary to determine the order of magnitude of $L_n$.
\subsubsection*{The case of a generic distribution $\rho$}

We will see in Section \ref{sec:distribution_generique}  that our proof  of Theorem \ref{th:encadrement} can in fact be adapted to any purely discrete distribution, not only distributions on $\mathbb{N}$. Finally, by combining this with the Ver\v{s}ik-Kerov theorem, we are able to derive the asymptotic for any real-valued distribution.

\begin{theorem}\label{th:distribution_generique}
Let $\rho$ be an arbitrary probability distribution on $\mathbb{R}$. Let $(\rho_1,\rho_2)$ be two measures on $(\R,\mathcal{B}(\R))$ such that $\rho=\rho_1+\rho_2$, with $\rho_1$  atomless and $\rho_2$ purely discrete (\emph{i.e.} there exists a countable set $S$ such that $\rho_2(\R\setminus S)=0)$). 
Let $(X_i)_{i\ge 1}$ be i.i.d. random variables with law $\rho$ and $L_n=L(X_1,\ldots,X_n)$. 
\begin{enumerate}[label=(\roman*)]
    \item {\bf (Purely discrete case)} If $\rho_1(\R)=0$ \emph{i.e.} $\rho=\rho_2$ then $$
\frac{1}{2}\leq \liminf_{n\to+\infty}\frac{L_n}{f_n}\leq \limsup_{n\to+\infty}\frac{L_n}{f_n}
\leq 1 \qquad \text{a.s.}
$$
where 
\begin{equation*}
    f_t= \inf_{x\geq 0}\left\{ x+\sum_{i\in S}\frac{t\rho(\{i\})}{t\rho(\{i\})+x}\right\};
\end{equation*}
    \item {\bf (Generic case)} For every distribution $\rho$
    $$\lim_{n\to \infty } \frac{L_n}{2\sqrt{n}}=\sqrt{\rho_1(\R)} \qquad a.s.$$
\end{enumerate}
\end{theorem}

\subsection{Context and related results}

As already mentioned, the Hammersley problem has been studied for several variants and with very different tools. In this article, we use the particle system approach (in particular couplings with a \emph{Hammersley process}) which goes back to \cite{Hamm} and was fruitfully exploited decades later, notably by \cite{AldousDiaconis,seppalainen1997increasing,CatorGroeneboom}. The method we are going to use to prove Theorem \ref{th:encadrement} is more precisely inspired  by the sources/sinks approach for the Hammersley process, this trick was introduced in \cite{CatorGroeneboom} (see also \cite{NousAlea,ciech2019order,gerin2024ulam}). It seems that there were no existing result for the Hammersley problem in the case of non-uniform discrete distributions.
 A particular difficulty posed by our problem is that, because of the inhomogeneity of $p_i$'s, we need non-asymptotic bounds on the Hammersley process.

Even if the methods are different we mention that there exist  very precise results for the problem of the longest \emph{weakly} increasing subsequence for uniform and non-uniform discrete distributions. Let $m\geq 1$ be fixed and let $(X_n)_{n\geq 1}$ be i.i.d. and uniform in $\{1,\dots, m\}$. Let $\LNS_n$ be the longest weakly increasing subsequence of $(X_1,\dots,X_n)$, it is easy to see that $\LNS_n$ will grow like $n/m$.
A consequence of a shape theorem for Young diagrams due to Kerov \cite[Sec.3.4,Th.2]{kerov2003asymptotic} 
(see also \cite[Th.4]{tracy2001distributions}) is that 
\begin{equation}\label{eq:TracyWidom}
\LNS_n= \frac{n}{m}+\sqrt{2n/m}\times \Lambda +\mathrm{o}_{\mathbb{P}}(\sqrt{n}),
\end{equation}
where $\Lambda$ is distributed like the largest eigenvalue in the Gaussian Unitary
Ensemble (see Johansson \cite[Th.1.7]{johansson2001discrete} and \cite{breton} for   different regimes where $m=m_n$ may depend on $n$). 
In the non-uniform case we have more generally that $\LNS_n$ grows like $n\times \max_i\{p_i\}$, it is proved in 
\cite[Sec.3]{houdre2009longest} (see also \cite{houdre2013limiting} for further improvements) that
\begin{equation}\label{eq:Houdre}
\LNS_n= n\times \max_i\{p_i\}+\sqrt{n}\times \mathcal{B}({\mathbf p}) +\mathrm{o}_{\mathbb{P}}(\sqrt{n}),
\end{equation}
where $\mathcal{B}({\mathbf p})$ is  a Brownian functional (its distribution strongly depends  on the number of maximal $p_i$'s).

Let $(P,Q)$ be the pair of random Young Tableaux associated to $(1,X_1),\dots,(n,X_n)$ through the RSK correspondence\footnote{For a review on the RSK correspondence we refer to Chap.7 in \cite{stanley1999enumerative}, notably Cor.7.23.11 and Th.7.23.17 regarding the relation between RSK and longest increasing subsequences.}.
The proofs of \cref{eq:TracyWidom,eq:Houdre} rely on the fact that $\LNS_n$ has the same distribution as the length of the first rows of $P$. It turns out that $\LIS_n$ has the same distribution as the number of rows in $P$, one may wonder if this alternative representation of $\LIS_n$ could provide additional information on its distribution. Given the wide variety of asymptotic behaviors, it seems to us that it would be very difficult to obtain a result similar to Theorem \ref{th:encadrement} with the RSK approach.

\bigskip

We  conclude this discussion with a word on fluctuations. The Baik–Deift–Johansson theorem states that
\begin{equation*}
\mathcal{L}_n=2\sqrt{n}+n^{1/6}\mathrm{TW}+o_{\mathbb{P}}(1),
\end{equation*}
where $\mathrm{TW}$ follows the Tracy-Widom distribution (see \cite[Sec.2]{Romik}). We leave the question of the order of magnitude of fluctuations of $L_n$ completely open, there is probably a wide variety of different asymptotics. However we mention that there is a general upper bound for the variance.

\begin{proposition}\label{prop:variance}
 Let $\rho$ be an arbitrary distribution, and let $X_1,\dots,X_n$ be i.i.d. random variables with common distribution $\rho$. Then
 $$
 \mathrm{Var}\left(L(X_1,\dots,X_n)\right)\le \E\left[L(X_1,\dots,X_n)\right].
 $$
\end{proposition}

\begin{proof}[Proof of Proposition \ref{prop:variance}]
The sequence of random variables $L_n$ is written as a sequence of \emph{self-bounding} functions, in the sense of \cite[Sec.3.3.]{boucheron2013concentration} (see also Example 3.11 therein). It follows then from the Efron-Stein inequality 
(see \cite[Cor.3.7]{boucheron2013concentration} for details) that for every $n\in\N$, $\mathrm{Var}(L_n)\le \E[L_n]$.
\end{proof}

\subsection{Description of the paper and strategy of proof}

\begin{itemize}
    \item In Section \ref{sec:Cv_ps_Theorem} we prove \cref{eq:Cv_ps_Theorem} in  Theorem \ref{th:encadrement}: the a.s. convergence of $L_n/\mathbb{E}[L_n]$.
 The proof is short and relies on classical tools from concentration measure theory applied to $L_n$. (Let us point out that inhomogeneity prevents the use of subadditivity, which is usually the essential tool for showing a.s. convergence.)
    \item Section \ref{sec:Hammersley} is devoted to the proof of Theorem \ref{th:poisson_domination}: we obtain  non-asymptotic bounds for the poissonized version of $L_n$. To this end, we introduce a dynamical representation of $L_n$ and obtain these bounds through couplings with an inhomogeneous version of discrete-time continuous-space  \emph{Hammersley process}. \\
    The construction of this process involves a free parameter $\alpha$. Optimizing over this $\alpha$ explains the variational formula of \cref{eq:fn}.
    \item In Section \ref{Sec:Properties_f_n} we collect some simple properties of $f_t$ and $w_t$ that will be used later.
    \item Section \ref{sec:PreuvePrincipale} is devoted to the end of the proof of Theorem \ref{th:encadrement}. The main part is to de-poissonize the bounds of Theorem \ref{th:poisson_domination}, using concentration inequalities.
    \item In Section \ref{sec:bornes_f} we obtain several estimates of $f_t$.
    \begin{itemize}
        \item In Section \ref{sec:greedy} we compare $f_t$ with the output of a naive greedy algorithm for finding an increasing subsequence. This analysis provides estimates which will be used in the proofs of Proposition \ref{prop:lighttail} and Proposition \ref{prop:heavytail}.
        \item In Section \ref{sec:bounds_ft} we obtain tight bounds for certain quantities related to $w_t,f_t$ with further assumptions on ${\mathbf p}$.
        \item In Section \ref{sec:explicit} we combine previous results to prove Proposition \ref{prop:lighttail} for geometric and Poisson distribution and Proposition \ref{prop:heavytail} for heavy-tailed distributions.
    \end{itemize}
    \item In Section \ref{sec:distribution_generique}, we show how to adapt the arguments of  the proof of Theorem \ref{th:encadrement} to the case of a generic distribution $\rho$.
\end{itemize}

\section{Preliminaries: convergence of $L_n/\mathbb{E}[L_n]$}\label{sec:Cv_ps_Theorem}

This short section is devoted to the proof of \cref{eq:Cv_ps_Theorem} in  Theorem \ref{th:encadrement}:  $L_n/\mathbb{E}[L_n]\to 1$ a.s..

\begin{proof}[Proof of \cref{eq:Cv_ps_Theorem}]
We first observe that if  {\bf p} has a finite support then  \cref{eq:Cv_ps_Theorem} is immediate.
Thus we assume from now on that the support of {\bf p} is infinite, in particular this implies that $\mathbb{E}[L_n]\to +\infty$. (This can be proved either by observing that $L_n$ is greater than the number of records or by applying Proposition \ref{prop:R_n_r_n} below, whose proof is independent from other sections.)

Using Proposition \ref{prop:variance} and Chebytchev's inequality we have for all $\eps>0$
    $$
\P\left(\left|\frac{L_n}{\E[L_n]}-1\right|\geq \eps\right)\le \frac{1}{\eps^2\E[L_n]}.
    $$
Since $\E[L_n]\to +\infty$ and  $L_{n+1}-L_n\le 1$, 
 there exists for every $n\geq 1$ an integer $m_n$ such that $n^2\le \E[L_{m_n}] \le n^2+1$. Using the  Borel-Cantelli Lemma
$\frac{L_{m_n}}{\E[L_{m_n}]}\to 1$  a.s.. Besides, for $m_n\le i \le m_{n+1}$
    $$
\frac{L_i}{\E[L_i]}\le \frac{L_{m_{n+1}}}{\E[L_{m_n}]}\le \frac{L_{m_{n+1}}}{\E[L_{m_{n+1}}]}\frac{(n+1)^2+1}{n^2} \stackrel{n\to+\infty}{\to} 1 \qquad \text{a.s.}.
    $$
    Similarly we have
$$
    \frac{L_i}{\E[L_i]}\geq \frac{L_{m_n}}{\E[L_{m_{n+1}}]}\geq \frac{L_{m_n}}{\E[L_{m_n}]}\frac{n^2}{(n+1)^2+1} \stackrel{n\to+\infty}{\to} 1 \qquad \text{a.s.},
    $$
    and \cref{eq:Cv_ps_Theorem} is proved.

\end{proof}

\begin{remark}\label{rem:Talagrand}
With a little care
one could get a more quantitative version of the convergence in \cref{eq:Cv_ps_Theorem} and obtain concentration inequalities for $L_n$ around its mean. Indeed $L_n$ satisfies the hypothesis of \cite[Th.7.1.2]{talagrand1995concentration} which yields concentration for $L_n$ around its median. Besides it is easy to check that $\mathrm{median}(L_n)/\mathbb{E}[L_n]\to 1$.
\end{remark}

\section{Poissonization and coupling with an inhomogeneous  Hammersley process}\label{sec:Hammersley}


The goal of this section is to prove Theorem \ref{th:poisson_domination}  which gives non-asymptotic bounds for a poissonized version of   $L_{n}$.

\begin{theorem}[Stochastic bounds for $L_{N_t}$]\label{th:poisson_domination} 
Let $\alpha, t >0$ be fixed. Let $(X_i)_{i\ge 1}$ be an i.i.d sequence with discrete distribution $\mathbf{p}$ and let $N_t$ be a Poisson random variable with mean $t$ independent of the sequence $(X_i)_{i\ge 1}$. We denote $\LIS_{N_t}$ the longest increasing subsequence of $(X_1,\ldots,X_{N_t})$. We can construct on the same probability space
\begin{itemize}
    \item two sequences of independent random variables  $(I^{(\alpha)}_i)_{i \geq 1}$, $(J^{(\alpha)}_i)_{i \geq 1}$ where  $I^{(\alpha)}_i\stackrel{\text{(d)}}{=}J^{(\alpha)}_i\stackrel{\text{(d)}}{=}\mathrm{Bernoulli}\left(\frac{p_i}{\alpha+p_i}\right)$;
    \item a random variable $\mathcal{H}^{(\alpha)}_t\stackrel{\text{(d)}}{=}\mathrm{Poisson}(\alpha t)$
\end{itemize}
such that the following holds:
\begin{enumerate}[label=(\roman*)]
    \item {\bf Upper bound.}  
    \begin{equation}\label{eq:Majo_E_LN}
    L_{N_t} \leq \mathcal{H}^{(\alpha)}_t+ \sum_{i\geq 1} I^{(\alpha)}_i.
    \end{equation}
    \label{item:poisson_upper}
    \item {\bf Lower bound.} 
     \begin{equation}\label{eq:Mino_E_LN}
     L_{N_t} \geq \min\set{\mathcal{H}^{(\alpha)}_t\ ;\ \sum_{i\geq 1} J^{(\alpha)}_i}.
    \end{equation}
     \label{item:poisson_lower}
\end{enumerate}
\end{theorem}

The proofs of \cref{eq:Majo_E_LN,eq:Mino_E_LN} relies on a coupling with an inhomogeneous variant of the Hammersley process (discrete time and continuous space). 
This process is parametrized by ${\bf p}$ and a parameter $\alpha>0$. Since this parameter $\alpha$ is free, we can optimize the RHS of \cref{eq:Majo_E_LN,eq:Mino_E_LN}, this explains the variational formulation of \cref{eq:fn,eq:wn}.

We have drawn inspiration from similar approaches for the Hammersley process, or some of its variants (see for instance \cite{CatorGroeneboom,AldousDiaconis,FerrariMartin,NousAlea,ciech2019order,gerin2024ulam}). 

\subsection{Preliminaries for the proof of Theorem \ref{th:poisson_domination}: coupling with a Hammersley process}

Fix $t>0$ throughout the section.
Let $\Pi_t$ be the random set $\Pi_t:=\bigcup_{i=1}^{+\infty} \Pi^i$ where $\Pi^i$'s are independent and each $\Pi^i$ is a homogeneous Poisson Point Process (PPP) with intensity $p_i$  on $[0,t]\times\{i\}$. 
For $i_0>1$ and $0\leq x\leq t$ we also put 
$$
\Pi^{i_0}_x=\Pi_t \cap \left( [0,x]\times\{i_0\}\right),\qquad 
\Pi^{\leq i_0}_x=\Pi_t \cap \left( [0,x]\times\{1,\dots ,i_0\}\right).
$$
For $\mathcal{P}$ a finite subset of $(\R_+)^2$, let define 
$$L(\mathcal{P}):=\max\{k\ge 0, \exists (t_i,x_i)_{1\le i\le k} \in \mathcal{P}^k, t_1<\ldots<t_k, x_1<\ldots<x_k\}$$
which is the maximal number of points of $\mathcal{P}$ an increasing function can get through.
By the superposition  property of PPP's we have
$$
L\left(\Pi_t\right) \stackrel{\text{(d)}}{=} \LIS(X_1,\dots,X_{N_t})\stackrel{\text{(d)}}{=} L_{N_t},
$$
where $N_t$ is a Poisson random variable with mean $t$.

\begin{figure}
\begin{center}
\includegraphics[width=150mm]{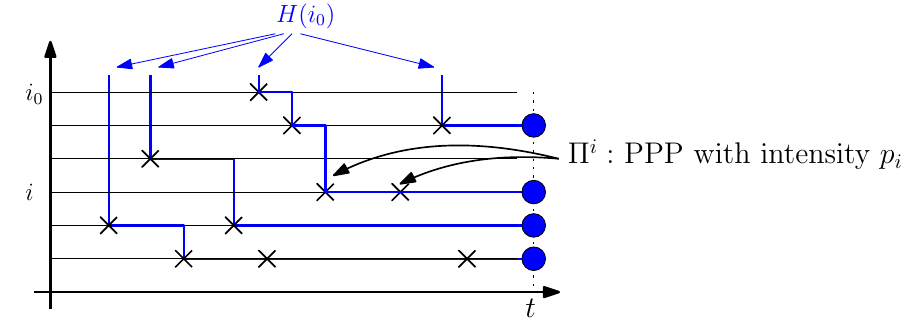} 
\end{center}
\caption{A sample of $\Pi_t$ and of the corresponding process $(H(i))_{i\geq 1}$. Observe that $L\left( \Pi^{\leq i_0}_{t}\right)=4=\mathrm{card}(H(i_0))$, as stated by Proposition \ref{prop:NombreLignes}.}
\label{fig:schemasHammersley}
\end{figure}

For every $i_0\in \{0,1,2,\dots \}$ the function $x\in[0,t] \mapsto L(\Pi^{\leq i_0}_x)$ is a non-decreasing integer-valued function whose all steps are equal to $+1$. Therefore this function is completely determined by the finite set
$$
H(i_0):=\set{x\leq t, L(\Pi^{\leq i_0}_x)=L(\Pi^{\leq i_0}_{x^-})+1}.
$$
(See \cref{fig:schemasHammersley}.) 
Sets $H(1),H(2),H(3),\dots,$ are finite subsets of $[0,t]$ whose elements are considered as particles. In the case where $p_i$'s are identical, it was observed by Hammersley \cite{Hamm} that the process $(H(i))_{i\geq 1}$ is a Markov process taking their values in the family of point processes of $[0,t]$.

For general $p_i$'s the process $(H(i))_{i\geq 1}$ is a inhomogeneous Markov process.
Similarly to the classical Hammersley process (see \cite[Sec.9]{Hamm} or \cite{AldousDiaconis}) the individual dynamic of particles is very easy to describe:
\begin{itemize}
\item {\bf The process $H$: particles evolution. (See \cref{fig:DynamiqueH}.)}  We put $H(0)=\emptyset$. In order to define $H(i+1)$ from $H(i)$ we consider particles of $H(i)$ from left to right. A particle at $x$ in $H(i)$ moves at time $i+1$ at the location of the leftmost available point $y$ in $\Pi^{i+1}_x$ (if any, otherwise it stays at $x$). This point $y$ is not available anymore for subsequent particles, as well as every other point of $\Pi^{i+1}_x$. \\
\begin{figure}
\begin{center}
\includegraphics[width=110mm]{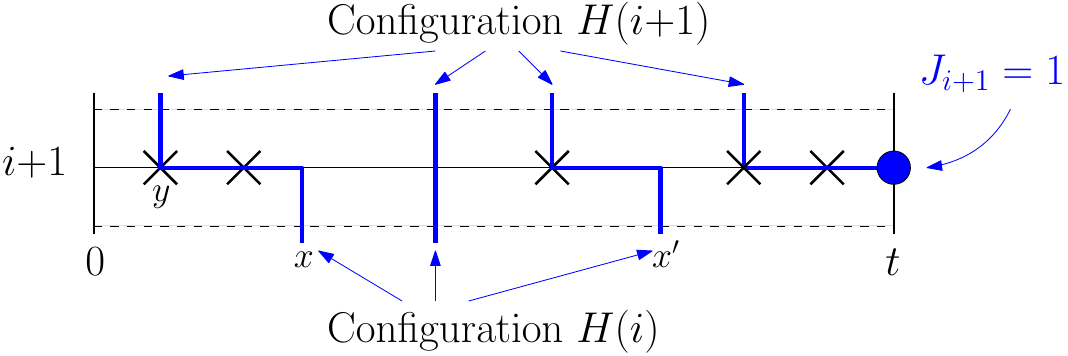} 
\end{center}
\caption{How to construct  $H(i+1)$ from  $H(i)$ and  $\Pi^{i+1}$ (points of $\Pi^{i+1}$ are depicted with $\times$'s as before). In this example $J_{i+1}=1$, which is indicated by the appearance of a new particle \textcolor{blue}{$\bullet$}.}
\label{fig:DynamiqueH}
\end{figure}
If there is a point in $\Pi^{i+1}$ which is on the right of $x':=  \max \{ H(i)\}$ then a new particle is created in $H(i+1)$, located at the leftmost point in $\Pi^{i+1}\cap (x',t)$ and we set $J_{i+1}=1$ (In our pictures this is indicated with a \textcolor{blue}{$\bullet$}.) If no particle is created , put $J_{i+1}=0$.
\end{itemize}

The process $(H(i))_i$ is designed in such a way that it keeps track of the length of longest increasing path in $\Pi_t$. 

\begin{proposition}[(\cite{Hamm}, Sec.9), see also (\cite{gerin2024ulam}, Prop.4)]\label{prop:NombreLignes}
For every $t \in (0,+\infty)$ and $i_0\geq 1$,
$$
L\left( \Pi^{\leq i_0}_{t}\right) = \mathrm{card}(H(i_0)).
$$
\end{proposition}

\subsection{Sources and sinks: stationarity}

In order to exploit Proposition \ref{prop:NombreLignes} we need to understand the asymptotic behaviour of processes $(H(i))_{i\geq 1}$. 

It follows from minor adjustments in previous results \cite{AldousDiaconis,Ferrari} that for every fixed parameter $\alpha>0$ the homogeneous PPP with intensity $\alpha$ on $\mathbb{R}$ is stationary for  $(H(i))_i$, as in the case of constant $p_i$'s.
Here we need non-asymptotic estimates for $(H(i))_i$.
 To solve this issue we use the trick of \emph{sources/sinks}  introduced formally and exploited by Cator and Groeneboom \cite{CatorGroeneboom} for the continuous Hammersley process:
\begin{itemize}
\item \emph{Sources} form a finite subset of $[0,t]\times \{0\}$ which plays the role of the initial configuration $H(0)$.  
\item \emph{Sinks} are points of $\{0\}\times  \N $  which add up to $\Pi$ when one defines the dynamics of $H(i)$.
\end{itemize}

\begin{figure}
\begin{center}
\includegraphics[width=150mm]{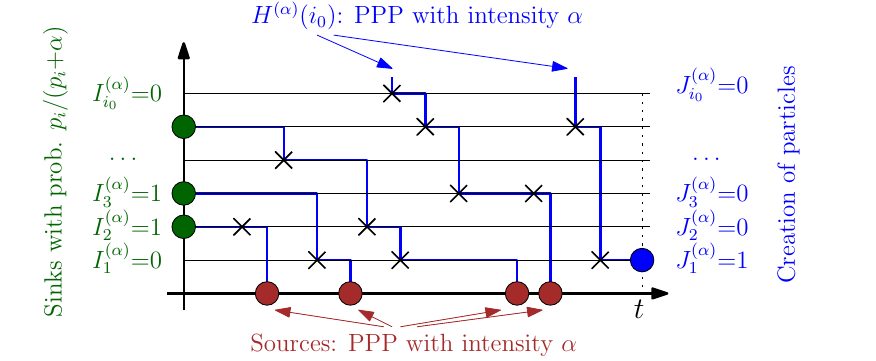} 
\end{center}
\caption{The same sample of $\Pi$ as \cref{fig:schemasHammersley}, with additional sources \textcolor{Brown}{$\bullet$} and sinks \textcolor{OliveGreen}{$\bullet$} ($I^{(\alpha)}_2=I^{(\alpha)}_3=I^{(\alpha)}_5=1$ while other $I^{(\alpha)}_i$'s are equal to zero.) 
In blue: the corresponding process $(H^{(\alpha)}(i))_{i\geq 1}$.\\
Proposition \ref{prop:Stationnaire<} states that 
(i) $H^{(\alpha)}(i_0)$ is also distributed as a homogeneous PPP with intensity $\alpha$\\ (ii) creations of particles on the right are distributed as sinks.
}
\label{fig:schemasHammersley_sources}
\end{figure}
\begin{figure}
\begin{center}
\includegraphics[width=140mm]{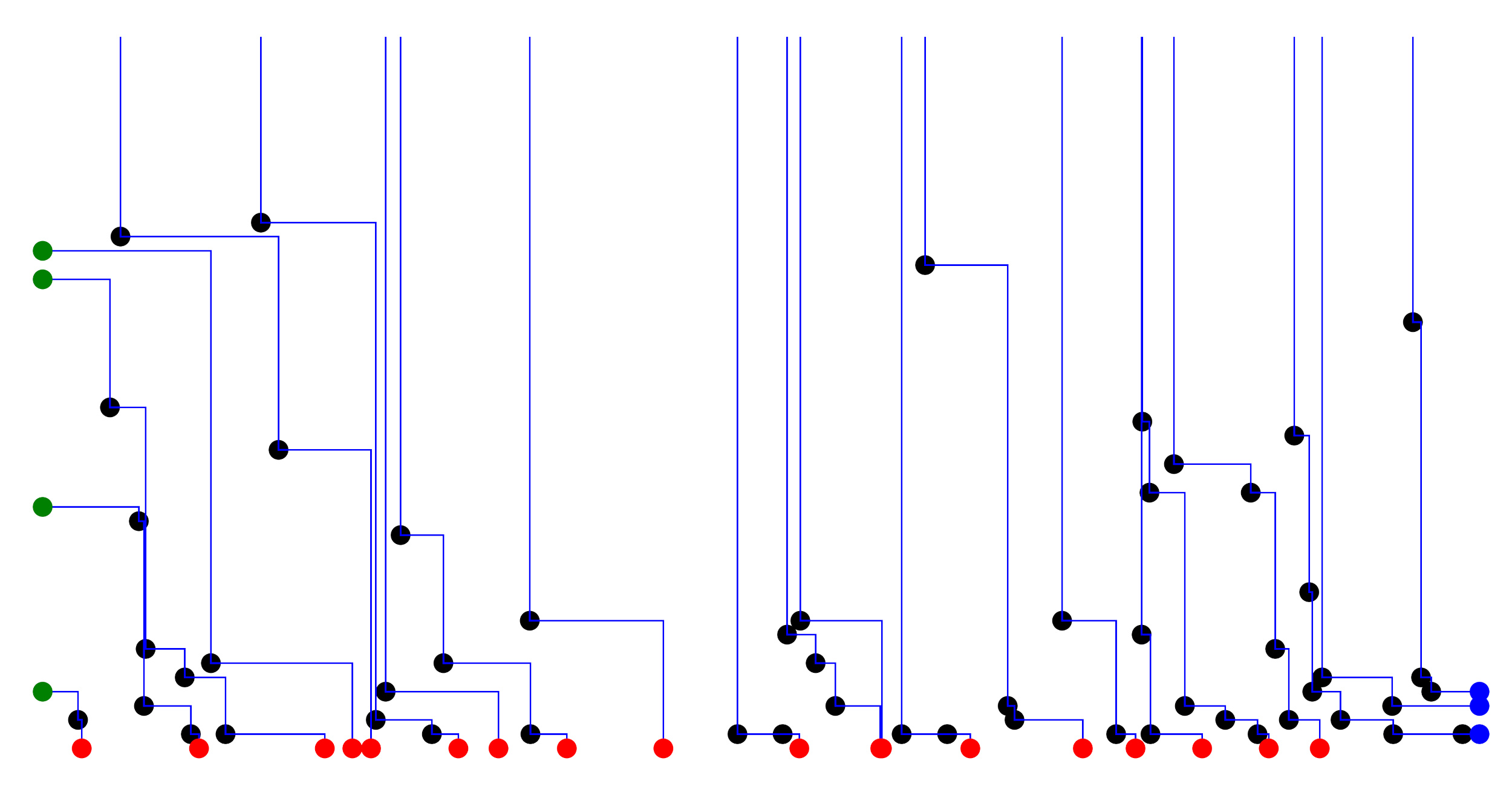}
\end{center}
\caption{A simulation of the inhomogeneous Hammersley process $(H^{(\alpha)}(i))_{i\geq 0}$ with a distribution $p_i\asymp i^{-1.2}$, $\alpha=0.2$, sources and sinks distributed as in Proposition \ref{prop:Stationnaire<}, $n=100$. Time goes from bottom to top, points of $\Pi$ are depicted with $\bullet$'s. As in \cref{fig:schemasHammersley_sources}, sources/sinks/creations of particles are respectively depicted with \textcolor{red}{$\bullet$}/\textcolor{ForestGreen}{$\bullet$}/\textcolor{blue}{$\bullet$}.\\
Proposition \ref{prop:Stationnaire<} states  that:\\
(i) Locations of particles at the top of the box are distributed as a homogeneous PPP with intensity~$\alpha$.\\
(ii) Spots \textcolor{blue}{$\bullet$} are distributed as \textcolor{ForestGreen}{$\bullet$}.}
\label{fig:SimuHammersley_alpha0-1}
\end{figure}

\begin{figure}
\begin{center}
\includegraphics[width=140mm]{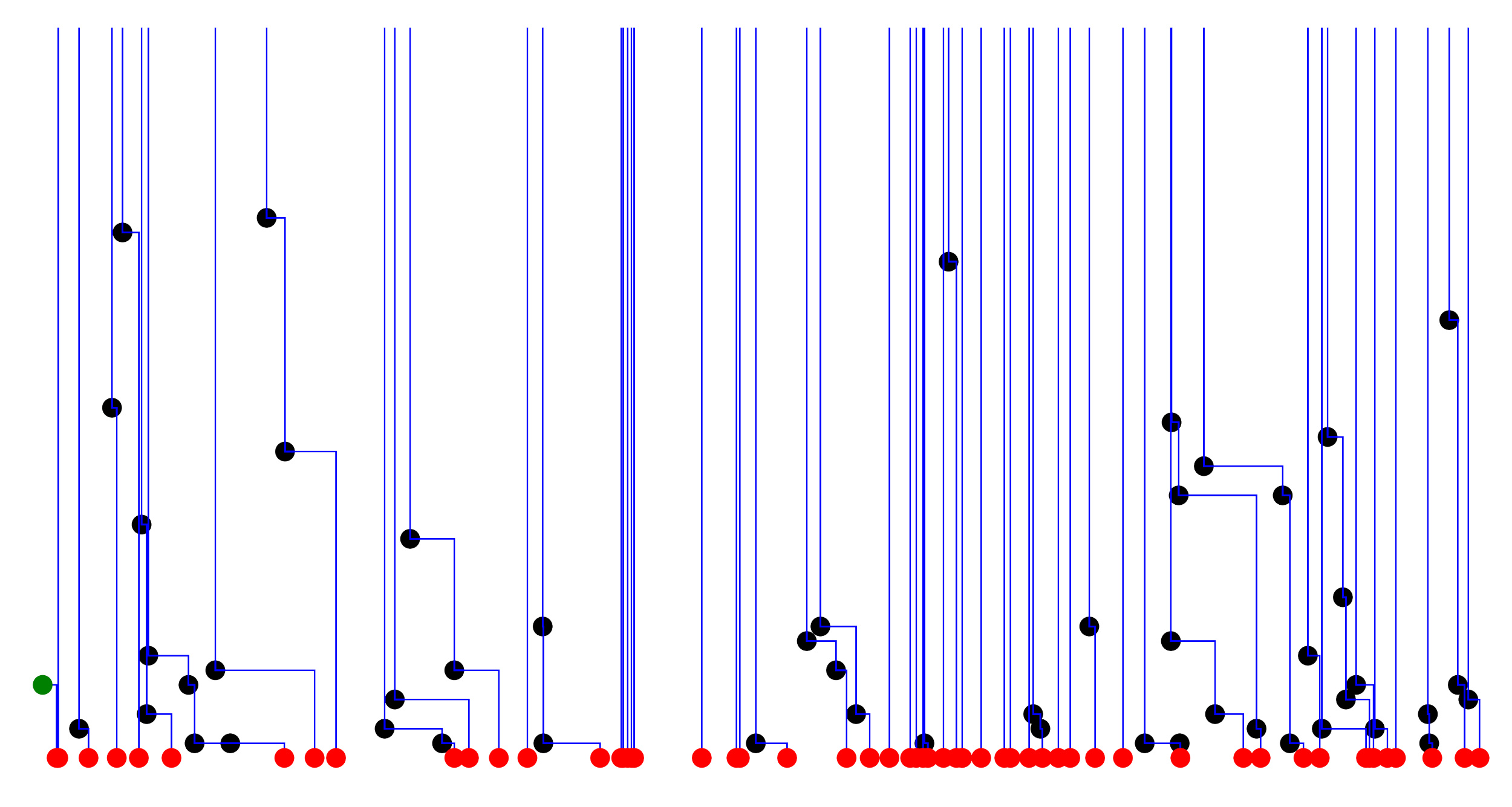}
\end{center}
\caption{A simulation of the inhomogeneous Hammersley process $(H^{(\alpha)}(i))_{i\geq 0}$ with the same realization  of $\Pi$ as in \cref{fig:SimuHammersley_alpha0-1}, with  $\alpha=0.6$.}
\label{fig:SimuHammersley_alpha0-6}
\end{figure}

Examples of dynamics of $H$ under the influence of sources/sinks is illustrated in  \cref{fig:schemasHammersley_sources}, see also  simulations in \cref{fig:SimuHammersley_alpha0-1,fig:SimuHammersley_alpha0-6}. Here is the inhomogeneous  discrete-time analogous of \cite[Th.3.1.]{CatorGroeneboom} (an important difference is that because of inhomogeneity of $p_i$'s we cannot expect to get a stationary version of the process $(H^{(\alpha)}(i))_{i\geq 0})$). 
\begin{proposition}[Steadiness of the Hammersley process with sources/sinks]\label{prop:Stationnaire<}
For every $t>0$, $\alpha>0$, let $(H^{(\alpha)}(i))_{i\geq 0}$ be the Hammersley process defined as $H$ with the following sets of sources and sinks.
\begin{itemize}
\item The set $H^{(\alpha)}(0)$ of sources is distributed according to a homogeneous PPP with intensity $\alpha$ on $[0,t]\times \{0\}$;
\item Let $I^{(\alpha)}_1,I^{(\alpha)}_2,\dots$ be a sequence of independent Bernoulli random variable where each $I_i$ has mean $p_i/(p_i+\alpha)$. Then we put a sink at position $(0,i)$ for each $i$ such that $I^{(\alpha)}_i=1$.
\item Sources, sinks, and $\Pi_n$ are independent. 
\end{itemize}
Then:
\begin{itemize}
\item[(i)]({\bf Bottom-to-top steadiness.}) For every $i$, the configuration
$H^{(\alpha)}(i)$ is a PPP with mean $\alpha$.
\item[(ii)]({\bf Left-to-right stationarity.}) Set $J^{(\alpha)}_i=1$ if a particle is created at time $i$, $J^{(\alpha)}_i=0$ otherwise. Then 
$$
(J^{(\alpha)}_1,J^{(\alpha)}_2,\dots ) \stackrel{\text{(d)}}{=} (I^{(\alpha)}_1,I^{(\alpha)}_2,\dots )
\stackrel{\text{(d)}}{=}
\mathrm{Bernoulli}(\tfrac{p_1}{p_1+\alpha})\otimes
\mathrm{Bernoulli}(\tfrac{p_2}{p_2+\alpha})\otimes \dots
$$
\end{itemize}
\end{proposition}

\begin{figure}
\begin{center}
\includegraphics[width=115mm]{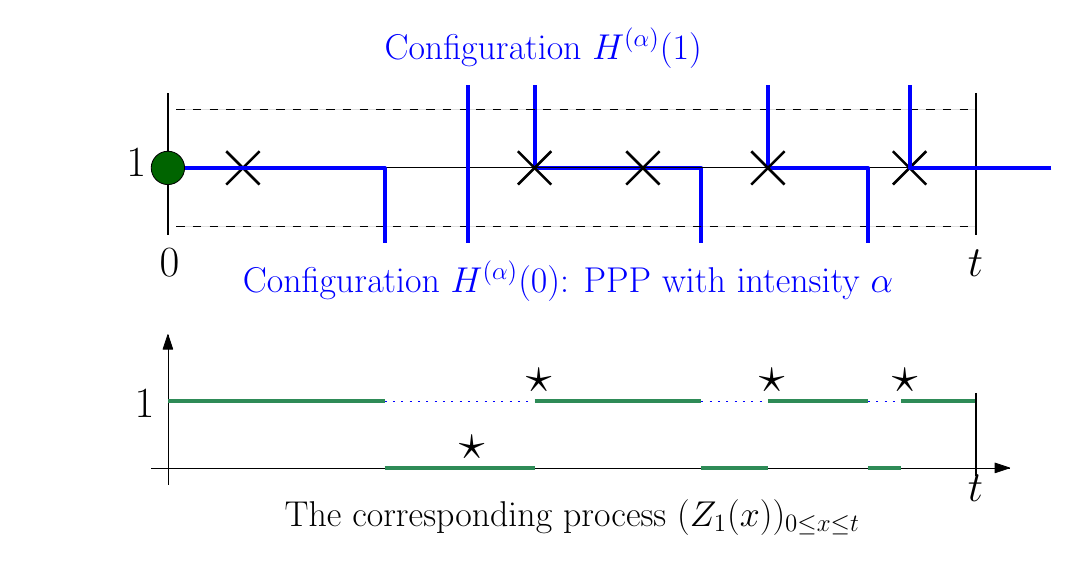}
\end{center}
\caption{Notation of Proposition \ref{prop:Stationnaire<}.}
\label{fig:ProcessH}
\end{figure}

\begin{proof}[Proof of Proposition \ref{prop:Stationnaire<}]

By induction it suffices to prove the proposition for $i=1$. 
Consider the process $(Z_1(x))_{0\leq x\leq t}$ defined by
$$
Z_1(x)=
\begin{cases}
 1&\text{if a Hammersley line touches }(x^+,1),\\
 0&\text{otherwise.}
 \end{cases}
$$

We regard $(Z_1(x))_{0\leq x\leq t}$ as a random process with time going from left to right. 
The initial value $Z_1(0)$ is equal to one if and only if there is a sink at $(0,1)$, which happens with probability $p_1/(p_1+\alpha)$. The process $(Z_1(x))_{0\leq x\leq t}$ is a continuous time Markov chain in $\set{0,1}$ with with '$+1$ rate' equal to $p_1$ and '$-1$ rate' equal to $\alpha$ (see a \cref{fig:ProcessH}). 
The Bernoulli$(p_1/(p_1+\alpha))$ distribution is stationary for this Markov chain, so that $(Z_1(x))_{0\leq x\leq t}$ is a stationary continuous-time Markov chain.

The configuration $H^{(\alpha)}(1)$  is given  by the union of (i) "$+1$ steps" of the process $Z_1$ and (ii) the points of $H^{(\alpha)}(0)$ that do not correspond to a '$-1$' jump (look at $\star$'s in \cref{fig:ProcessH}). It remains to prove that this set is distributed as a homogeneous PPP with intensity $\alpha$. This will follow from the following Lemma.

\begin{figure}
\begin{center}
\includegraphics[width=130mm]{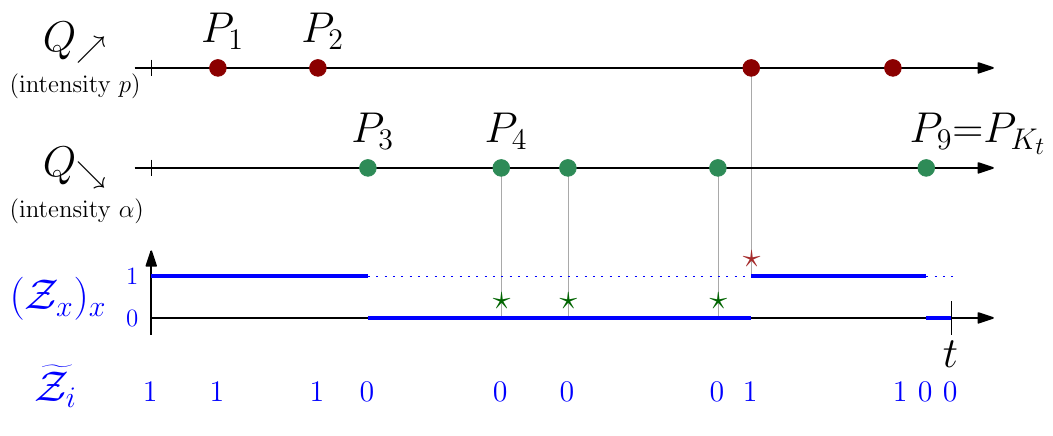} 
\end{center}
\caption{Notation of Lemma \ref{lem:Burke}. Points of $Q=Q_0\cup Q_1$ are depicted with $\star$'s. (Green $\star$'s correspond to $Q_0$, the brown $\star$ corresponds to $Q_1$.)}
\label{fig:H}
\end{figure}

\begin{lemma}\label{lem:Burke}
 (The reader is invited to look at Fig.\ref{fig:H} for notation.)\\
Let $t> 0$. Let $Q_\nearrow,Q_\searrow$ be two independent homogeneous PPP over $(0,t)$ with respective intensities $p,\alpha>0$. Let $(\mathcal{Z}_x)_{0\leq x\leq t}$ be the continuous-time Markov process with values in $\set{0,1}$ such that
\begin{itemize}
    \item $\mathcal{Z}_0$ is  independent from $Q_\nearrow,Q_\searrow$ and drawn accordingly to a Bernoulli$(p/(p+\alpha))$;
    \item If $\mathcal{Z}_{x^-}=0$ and $x\in Q_\nearrow$ then $\mathcal{Z}_{x}=1$;
    \item If $\mathcal{Z}_{x^-}=1$ and $x\in Q_\searrow$ then $\mathcal{Z}_{x}=0$.
\end{itemize}
(In words  '+1' jumps are given by $Q_\nearrow$ and its '-1' jumps are given by $Q_\searrow$.)

Let $Q_0,Q_1$ be respective point processes given by \emph{unused} "$-1$ jumps" and \emph{used} "$+1$ jumps":
$$
Q_0=\set{x\in Q_\searrow\text{ such that } \mathcal{Z}_{x^-}=0},\qquad 
Q_1=\set{x\in Q_\nearrow\text{ such that } \mathcal{Z}_{x^-}=0}.
$$
Then:
\begin{itemize}
    \item[(i)] the process $Q:=Q_0\cup Q_1$ is a homogeneous PPP with intensity $\alpha$;
    \item[(ii)] $\mathcal{Z}_{t}\in\set{0,1}$ is a $\mathrm{Bernoulli}(p/(p+\alpha))$ independent from the point process $Q$.
\end{itemize}
\end{lemma}

\begin{proof}[Proof of Lemma \ref{lem:Burke}]
The point process $Q_\nearrow\cup Q_\searrow$ is a homogeneous PPP with intensity $p+\alpha$, independent from $\mathcal{Z}_0$. We claim that $Q$ is a subset of $Q_\nearrow\cup Q_\searrow$ where each point in  $Q_\nearrow\cup Q_\searrow$ is taken independently with probability $\alpha/(p+\alpha)$, it is therefore a homogeneous  PPP with intensity $\alpha$, as required. 

We need a few notation in order to prove the claim. Set $P_0=0$ and for $i\geq 1$ let $P_i$ be the $i$-th leftmost point of $Q_\nearrow\cup Q_\searrow$ and let $(\widetilde{\mathcal{Z}}_i)_{i\geq 0}$ be the discrete-time embedded chain associated to $\mathcal{Z}$, \emph{i.e.} $\widetilde{\mathcal{Z}}_i=\mathcal{Z}_{P_i}$ for every $i$.
First observe that $\widetilde{\mathcal{Z}}_i$ has probability transitions given by $\begin{pmatrix} \frac{\alpha}{p+\alpha} & \frac{p}{p+\alpha}\\ \frac{\alpha}{p+\alpha} & \frac{p}{p+\alpha}\end{pmatrix}$, and is therefore a sequence of i.i.d. $\mathrm{Bernoulli}(p/(p+\alpha))$.

Now, for every $i\geq 1$,
\begin{align*}
\set{P_i\in Q}&=\set{P_i\text{ is a "used" $+1$ jump}}\cup \set{P_i\text{ is a "unused" $-1$ jump}}\\
&=\set{P_i\in Q_\nearrow, \widetilde{\mathcal{Z}}_{i-1}=0}\cup \set{P_i\in Q_\searrow, \widetilde{\mathcal{Z}}_{i-1}=0}\\
&=\set{\widetilde{\mathcal{Z}}_{i-1}=0},
\end{align*}
and therefore the sequence  $\set{P_i\in Q}$ is a sequence of i.i.d. $\mathrm{Bernoulli}_{\geq 0}(p/(p+\alpha))$, as claimed. It implies assertion (i): $Q$ is a homogeneous PPP with intensity $\alpha$.

Let now $P_{K_t}$ be the rightmost point of $Q_\nearrow\cup Q_\searrow$. The above construction also shows that no matter $Q_0,Q_1$, $\mathcal{Z}_t=(\widetilde{\mathcal{Z}}_{K_t})$ is a fresh Bernoulli$(p/(p+\alpha))$ which shows (ii).

\end{proof}
To conclude the proof of Proposition \ref{prop:Stationnaire<} we apply successively Lemma  \ref{lem:Burke} with
$$
\begin{array}{l l l}
p=p_1, & Q_\searrow=H^{(\alpha)}(0)=\set{\text{sources of $(0,n)$}},&  Q_\nearrow=\Pi_1,\\
p=p_2, & Q_\searrow=H^{(\alpha)}(1),&  Q_\nearrow=\Pi_2,\\
p=p_3, & Q_\searrow=H^{(\alpha)}(2),&  Q_\nearrow=\Pi_3,\\
\dots & & 
\end{array}
$$
\end{proof}

With Proposition \ref{prop:Stationnaire<} we are now ready to prove  the stochastic bounds \cref{eq:Majo_E_LN,eq:Mino_E_LN} of Theorem \ref{th:poisson_domination}.

\begin{proof}[Proof of Theorem \ref{th:poisson_domination}.]

For $x\leq t$, let ${\sf So}_x=H(0)\cap[0,x]$  be the random set of sources with abscissas $\leq x$  and for $i_0\geq 1$, let  ${\sf Si}_{i_0}$ be the set of sinks with ordinates $\leq i_0$. In particular,
$$
\mathrm{card}({\sf So}_x)\stackrel{\text{(d)}}{=}\mathrm{Poisson}(x\alpha ),\qquad \mathrm{card}({\sf Si}_{i_0})=\sum_{i\leq i_0} I^{(\alpha)}_i\stackrel{\text{(d)}}{=}\sum_{i\leq i_0}\mathrm{Bernoulli}\left(\frac{p_i}{p_i+\alpha}\right).
$$
 It is convenient to use the notation $L_{=<}(\mathcal{P})$
which is, as before, the length of the longest increasing path taking points in $\mathcal{P}$ but when the path is also allowed to go through several sources (which have however the same ordinate) or several sinks (which have the same abscissa). Formally,
$$
L_{=<}(\mathcal{P})=\max\set{L; \hbox{ there exists }P_1  =\prec P_2 =\prec \dots =\prec P_L, \text{ where each }P_i \in \mathcal{P} },\\
$$
where
$$
(a,b)=\prec (a',b') \text{ if }
\begin{cases}
 &a<a'\text{ and }b<b', \\
\text{ or } &a=a'=0 \text{ and }b<b', \\
\text{ or } &a<a' \text{ and }b=b'=0.
\end{cases}
$$
Proposition \ref{prop:NombreLignes} generalizes easily to the settings of sinks and sources.
\begin{claim}(See \cref{fig:schemasHammersley_sources}.)
\begin{align}
L_{=<}\left( \Pi^{\leq i_0}_t\cup {\sf So}_t\cup {\sf Si}_{i_0}\right) &=\mathrm{card}(H^{(\alpha)}(i_0))+\sum_{i\leq i_0} I^{(\alpha)}_i\label{claim:GaucheHaut}\\
&=\mathrm{card}(H^{(\alpha)}(0))+\sum_{i\leq i_0} J^{(\alpha)}_i.\label{claim:DroiteBas}
\end{align}
\end{claim}
(For the example of \cref{fig:schemasHammersley_sources} one has indeed that
$L_{=<}\left( \Pi^{\leq i_0}_t\cup {\sf So}_t\cup {\sf Si}_{i_0}\right)$ is equal to $5=2+3=$  $\mathrm{card}(H^{(\alpha)}(i_0))+\sum_{i\leq i_0} I^{(\alpha)}_i$  but also equal to $4+1= \mathrm{card}(H^{(\alpha)}(0)) +\sum_{i\leq i_0} J^{(\alpha)}_i$.) 
\begin{proof}[Proof of the Claim.] By the same reasoning as in the proof of Proposition \ref{prop:NombreLignes} the LHS is exactly the number of broken lines in the box $[0,t]\times [0,i_0]$. 
The number of such broken lines can be double-counted:
\begin{itemize}
    \item either by saying that each broken lines escapes the box either through the left (it thus corresponds to a sink) or through the top (and is thus counted by $\mathrm{card}(H^{(\alpha)}(i_0))$), hence $$L_{=<}\left( \Pi^{\leq i_0}_t\cup {\sf So}_t\cup {\sf Si}_{i_0}\right) =\mathrm{card}(H^{(\alpha)}(i_0))+\sum_{i\leq i_0} I^{(\alpha)}_i;$$
    \item or by saying that each broken lines enters the box either through the bottom (it thus corresponds to a source) or through the right (and is thus counted by the number of creations of particles) hence 
    $$
    L_{=<}\left( \Pi^{\leq i_0}_t\cup {\sf So}_t\cup {\sf Si}_{i_0}\right) =\mathrm{card}(H^{(\alpha)}(0))+\sum_{i\leq i_0} J^{(\alpha)}_i.
    $$
\end{itemize}
\end{proof}

\noindent{\bf Proof of \cref{eq:Majo_E_LN}: upper bound in Theorem \ref{th:poisson_domination}.}
Adding sources and sinks may not decrease longest increasing paths. Thus,
\begin{align}
L\left( \Pi^{\leq i_0}_{t}\right) &\leq L_{=<}\left( \Pi^{\leq i_0}_t\cup {\sf So}_t\cup {\sf Si}_{i_0}\right) \notag\\
&= \mathrm{card}(H^{(\alpha)}(i_0))+\sum_{i\leq i_0} I^{(\alpha)}_i\qquad \text{ (using \cref{claim:GaucheHaut})}.\label{eq:presque_Majo_E_LN}
\end{align}
Let now $i'$ be the ordinate of the highest point of $\Pi_t\cup {\sf So}_t\cup {\sf Si}_{i_0}$. Then by construction,  for every $i_0\geq i'$ one has $H^{(\alpha)}(i_0)=H^{(\alpha)}(i')$, in particular $H^{(\alpha)}(i_0)$ converges almost surely when $i_0\to +\infty$. Besides we have that, using bottom-to-top steadiness (Proposition \ref{prop:Stationnaire<} (i)),
$$
H^{(\alpha)}(i_0)\stackrel{\text{(d)}}{=}H^{(\alpha)}(0)\stackrel{\text{(d)}}{=}
\text{PPP with intensity $\alpha$ on $(0,t)$}.
$$
It follows that $H^{(\alpha)}(i')\stackrel{\text{(d)}}{=}
\text{PPP with intensity  $\alpha$ on $(0,t)$}$. 
Letting now $i_0\to +\infty$ in \cref{eq:presque_Majo_E_LN} we obtain
$$
L\left( \Pi_{t}\right) \leq \mathrm{card}(H^{(\alpha)}(i'))+ \sum_{i\geq 1} I^{(\alpha)}_i =: \mathcal{H}^{(\alpha)}_t+ \sum_{i\geq 1} I^{(\alpha)}_i,
$$
where $\mathcal{H}^{(\alpha)}_t$ is a Poisson random variable with mean $\alpha t$.

\noindent{\bf Proof of \cref{eq:Mino_E_LN}: lower bound in Theorem \ref{th:poisson_domination}.}

A maximizing path for $L_{=<}\left( \Pi^{\leq i_0}_t\cup {\sf So}_t\cup {\sf Si}_{i_0}\right)$ goes through some points of $\Pi^{\leq i_0}_t$ (at most $L( \Pi^{\leq i_0}_{t})$ of them) and goes through sinks or sources but not both. Hence 
\begin{align*}
L_{=<}\left( \Pi^{\leq i_0}_t\cup {\sf So}_t\cup {\sf Si}_{i_0}\right)&\leq L\left( \Pi^{\leq i_0}_{t}\right) +\max\set{ \mathrm{card}(H^{(\alpha)}(0)); \sum_{i=1}^{i_0}I^{(\alpha)}_i}.
\end{align*}
Assume that $\mathrm{card}(H^{(\alpha)}(0))\ge \sum_{i=1}^{i_0}I^{(\alpha)}_i$. Then we get
\begin{equation*}
    L\left( \Pi^{\leq i_0}_{t}\right) \ge L_{=<}\left( \Pi^{\leq i_0}_t\cup {\sf So}_t\cup {\sf Si}_{i_0}\right)- \mathrm{card}(H^{(\alpha)}(0)) \stackrel{\text{\cref{claim:DroiteBas}}}{=} \sum_{i=1}^{i_0}J^{(\alpha)}_i.
\end{equation*}
Assume on the contrary that $\mathrm{card}(H^{(\alpha)}(0))\leq \sum_{i=1}^{i_0}I^{(\alpha)}_i$. Using now \cref{claim:GaucheHaut}, $L\left( \Pi^{\leq i_0}_{t}\right)\ge  \mathrm{card}(H^{(\alpha)}(i_0))$. Finally
\begin{align*}
 L\left( \Pi^{\leq i_0}_{t}\right)\ge  \min\set{ \mathrm{card}(H^{(\alpha)}(i_0)); \sum_{i=1}^{i_0}J^{(\alpha)}_i}.
\end{align*}
As for the proof of the upper bound, we let $i_0\to +\infty$ and obtain
$$
L(\Pi_t) \geq \min\set{\mathcal{H}^{(\alpha)}_t,\sum_{i\geq 1}J^{(\alpha)}_i}.
$$
\end{proof}

\section{Basic properties of $f_t$ and $w_t$}\label{Sec:Properties_f_n}
We collect here for convenience some simple properties of both quantities $f_t$ and $w_t$ that will be used throughout Section \ref{sec:PreuvePrincipale}.
To this aim we put
\begin{equation}\label{eq_def_g}
    g_t(a)= at+\sum_{i\geq 1}\frac{p_i}{p_i+a}
\end{equation}
and recall the notation:
\begin{align*}
    f_t&= \inf_{x\geq 0}\left\{ x+\sum_{i\geq 1}\frac{tp_i}{tp_i+x}\right\}=
    \inf_{a\geq 0}\left\{ at+\sum_{i\geq 1}\frac{p_i}{p_i+a}\right\}=
    \inf_{a\geq 0}\left\{ g_t(a)\right\},\\
w_t&\text{ is the unique positive solution of }w_t=\sum_{i\geq 1}\frac{tp_i}{tp_i+w_t}.
\end{align*}

\begin{lemma}[Properties of $f_t,w_t$]\label{lem:basic}
For every $t>0$ and $\varepsilon\in (0,1)$,
    \begin{enumerate}[label=(\roman*)]
        \item $f_t \leq 2\sqrt{t}$ and $f_t=\mathrm{o}(\sqrt{t})$. \label{majorf}
         \item  $w_t\leq f_t\leq 2w_t$.\label{encadref}
       \item  $f_{t(1+\eps)}\leq (1+\eps)f_t.$\label{ineqf}
        \item  $(1-\eps) w_t \leq w_{t(1-\eps)}.$\label{ineqw}
        \item\label{item:majo_somme_degueu} We have the upper bound
$$  \sum_{i=1}^\infty\frac{tp_i}{(tp_i+w_t)^2}\leq 4.$$
    \end{enumerate}
\end{lemma}
\begin{proof}[Proof of Lemma  \ref{lem:basic}]
For \ref{majorf}, let us  write
    $$
f_t\leq g_t(\frac{1}{\sqrt{t}})=\sqrt{t}+\sum_{i=1}^{+\infty}\frac{p_i}{p_i+\frac{1}{\sqrt{t}}}\leq \sqrt{t}+\sqrt{t}\sum_{i=1}^{+\infty}p_i=2\sqrt{t}.
    $$
Moreover, for every $\eps>0$    
    $$
f_t\leq g_t(\frac{\varepsilon}{\sqrt{t}})=\varepsilon\sqrt{t}+\sum_{i=1}^{\varepsilon\sqrt{t}}\frac{p_i}{p_i+\frac{\varepsilon}{\sqrt{t}}}+\sum_{i> \varepsilon\sqrt{t}}\frac{p_i}{p_i+\frac{\varepsilon}{\sqrt{t}}}\leq\varepsilon\sqrt{t}+\varepsilon\sqrt{t}+\frac{\sqrt{t}}{\varepsilon}\sum_{i> \varepsilon\sqrt{t}}p_i\le 3\varepsilon\sqrt{t},
    $$    
for large enough $t$. Hence \ref{majorf}.

We turn to the proof of \ref{encadref}.
     We have that 
     $$
 f_t\leq g_t(\frac{w_t}{t})=w_t +    \sum_{i\geq 1}\frac{p_i}{p_i+\frac{w_t}{t}}=2w_t.
     $$
 Regarding the left inequality in \ref{encadref}, first observe that,
 since for every $t$, $g_t(\cdot)$ is strictly convex on $(0,\infty)$ with $\lim_{a\to 0^+} g_t(a)=\lim_{a\to +\infty} g_t(a)=+\infty$, there exists a unique $\alpha^\star_t\in (0,\infty)$ such that 
 \begin{equation}\label{eq:def_alpha_n_star}
 f_t=g_t(\alpha_t^\star).
 \end{equation}
Assume first  that $t\alpha_t^\star>w_t$, then 
     $f_t=t\alpha_t^\star+\sum_{i=1}^{+\infty}\frac{p_i}{p_i+\alpha_t^\star}> w_t$. Conversely, if $t\alpha_t^\star\leq w_t$, then
     $$
f_t\geq \sum_{i=1}^{+\infty}\frac{p_i}{p_i+\alpha_t^\star}\geq \sum_{i=1}^{+\infty}\frac{p_i}{p_i+\frac{w_t}{t}}=w_t.
     $$
Hence, in any case, $f_t\ge w_t$.

   In order to prove \ref{ineqf} we write
   $$
f_{t(1+\eps)}= \inf_{a>0}\left\{t(1+\eps)a+\sum_{i=1}^{+\infty}\frac{p_i}{p_i+a}\right\}\leq \inf_{a>0}\left\{t(1+\eps)a+(1+\eps)\sum_{i=1}^{+\infty}\frac{p_i}{p_i+a}\right\}= (1+\eps)f_t.
    $$
    
     To prove \ref{ineqw} we write for $t>0$
    $$
\frac{w_t}{t}=\frac{1}{t}\sum_{i=1}^\infty \frac{p_i}{p_i+\frac{w_t}{t}}.
    $$
    From this equality we see that $t\mapsto \frac{w_t}{t}$ is non-increasing. Hence we have for all $t> 0$ and $\eps\in (0,1)$
    $$
\frac{w_t}{t}\leq \frac{w_{(1-\eps)t}}{(1-\eps)t},
    $$
    and the inequality follows. 

Let us finally prove \ref{item:majo_somme_degueu}. 
Let $\alpha_t^\star$ be defined by \cref{eq:def_alpha_n_star}, we first observe that $\alpha_t^\star t\leq  g_t(\alpha_t^\star) =  f_t\le 2w_t$ (using  \ref{encadref}).
Besides $g_t'(\alpha_t^\star)=0$, \emph{i.e.} 
$$\sum_{i=1}^\infty\frac{tp_i}{(tp_i+t\alpha_t^\star)^2}=1.$$
Therefore
$$1\ge \sum_{i=1}^\infty\frac{tp_i}{(tp_i+2w_t)^2}\ge \frac{1}{4}\sum_{i=1}^\infty\frac{tp_i}{(tp_i+w_t)^2}.$$
\end{proof}

\section{Proof of Theorem \ref{th:encadrement}}\label{sec:PreuvePrincipale}
\subsection{Proof of  Theorem \ref{th:encadrement}: upper bound for $\mathbb{E}[\LIS_n]$}
\label{SectionUpper}

The aim of this section is to prove the upper bound for $\LIS_n$ given in Theorem \ref{th:encadrement}, namely
\begin{equation}\label{eq:upper_bound_Esperance}
\limsup_{n\to \infty} \frac{\mathbb{E}[\LIS_n]}{f_n}\leq 1.
\end{equation}

\begin{proof}[Proof of \cref{eq:upper_bound_Esperance}: upper bound in Theorem \ref{th:encadrement}]
The upper bound for the expectation of $\LIS_n$ is a direct consequence of the stochastic domination of $\LIS_{N_t}$ given in Theorem \ref{th:poisson_domination}  \ref{item:poisson_upper}. Indeed, let $\alpha_t^\star$ be such that 
$$f_t=g_t(\alpha_t^\star)=\alpha_t^\star t +\sum_{i\ge 1}\frac{p_i}{p_i+\alpha_t^\star}.$$
Choosing $\alpha=\alpha_t^\star$ in  Theorem \ref{th:poisson_domination} \ref{item:poisson_upper} implies that, for any $t> 0$, 
\begin{equation}\label{eq:DominationEsperance}
\E(\LIS_{N_t})\le f_t.
\end{equation}
Fix now $\eps>0$. By monotonicity of $\LIS_n$ it holds that
\begin{equation}\label{eq:OnUtiliseIndependance}
\E[\LIS_{N_{n(1+\eps)}}]\geq \E[\LIS_{N_{n(1+\eps)}}\mathbf{1}_{N_{n(1+\eps)}\geq n}]\geq \E[\LIS_n]\times\mathbb{P}(N_{n(1+\eps)}\geq n).
\end{equation}
Since $\mathbb{P}(N_{n(1+\eps)}\geq n)\stackrel{n\to+\infty}{\to} 1$, we get that, if  $n$ is large enough, 
$\E[\LIS_n]\leq (1+\eps)\E[\LIS_{N_{n(1+\eps)}}]$. We get
    $$
\frac{\E[\LIS_n]}{f_{n(1+\eps)}}\leq (1+\eps)\frac{\E[\LIS_{N_{n(1+\eps)}}]}{f_{n(1+\eps)}}\stackrel{\text{eq.}\eqref{eq:DominationEsperance}}{\leq} (1+\eps).
    $$
Using that 
$
f_{n(1+\eps)}\leq (1+\eps)f_{n}
$ (\emph{c.f.} Lemma \ref{lem:basic} \ref{ineqf}),
we get, for all $\eps$, if $n$ is large enough,
$$
\frac{\E[\LIS_n]}{f_n}\leq (1+\eps)^2
$$
which gives the upper bound of Theorem \ref{th:encadrement}.
\end{proof}

\subsection{
Proof of Theorem \ref{th:encadrement}: lower bound for $\mathbb{E}[\LIS_n]$ }
\label{SectionLower}

Recall that $w_t$ the unique positive solution of the equation  
$$w_t=\sum_{i=1}^\infty \frac{tp_i}{tp_i+w_t}.$$
The aim of this section is to prove the lower bound in Theorem \ref{th:encadrement}, namely
\begin{equation}\label{eq:lower_bound_Theorem}
\liminf_{n\to \infty} \frac{\mathbb{E}[\LIS_n]}{w_n}\geq 1.
\end{equation}
We first prove a lower bound for the poissonized version $L_{N_t}$.

\begin{proposition}\label{prop:concentration2}
Let $N_t$ be a $\mathrm{Poisson}$ random variable with mean $t$.
   For all $\eps\in (0,1)$,  we have:
    $$
\mathbb{P}(L_{N_t}<(1-\eps)w_t)\leq 2\exp\left(-\frac{\varepsilon^2 w_t}{40}\right).
    $$
\end{proposition}

\begin{proof}
\newcommand{\alfa}{\widetilde{\alpha}_t}

Using Theorem \ref{th:poisson_domination} \ref{item:poisson_lower}  with  $\alfa=w_t/t$, we get that
 $$L_{N_t} \geq \min\set{\mathcal{H}^{(\alfa)}_t\ ;\ \sum_{i\geq 1} J^{(\alfa)}_i}$$
 where    $(J^{(\alfa)}_i)_{i \geq 1}$ are independent random variables with Bernoulli distributions with respective means 
 $\left(\frac{p_i}{\alfa+p_i}\right)=\left(\frac{t p_i}{w_t+t p_i}\right)$ and 
    $\mathcal{H}^{(\alfa)}_t$ is a Poisson random variable with parameter $t\alfa =w_t$.
Thus
\begin{equation*}
    \mathbb{P}(L_{N_t}<(1-\eps)w_t)\leq \mathbb{P}(\mathcal{H}^{(\alfa)}_t\le (1-\eps)w_t)+ \mathbb{P}(\sum_{i\ge 1} J^{(\alfa)}_i\le (1-\eps)w_t).
\end{equation*}
From a standard concentration inequality for Poisson random variables \cite[eq.(2.5) and Rem.2.6]{janson2011random} we have
\begin{equation}\label{eq:lowerbound1}
\mathbb{P}(\mathcal{H}^{(\alfa)}_t\le (1-\eps)w_n)\le \exp\left(-\frac{\varepsilon^2 w_t}{2}\right).
\end{equation}
Let $Z_t:=\sum_{i\ge 1} J^{(\alfa)}_i$ and for $M>0$,  $Z^M_{t}=\sum_{i=1}^M J^{(\alfa)}_i$.
Let choose $M:=M(t)$ such that 
$$(1-\varepsilon) w_t-\E(Z^M_{t})=\sum_{i>M}\frac{tp_i}{tp_i+w_t}-\varepsilon w_t\le -\frac{\varepsilon w_t}{2}.$$
We get
$$\mathbb{P}(Z_{t}\le (1-\eps)w_t)\le \mathbb{P}(Z^M_{t}\le (1-\eps)w_t)\le\mathbb{P}\left(Z^M_{t}-\E(Z^M_{t})\le -\frac{\varepsilon w_t}{2}\right).$$
Using Bernstein inequality, we get 
$$\mathbb{P}(Z^M_{t}\le (1-\eps)w_t)\le \exp\left(-\frac{\frac{(\varepsilon w_t)^2}{8}}{w_t\sum_{i=1}^M\frac{tp_i}{(tp_i+w_t)^2}+\frac{\varepsilon w_t}{3}}\right)\le \exp\left(-\frac{\frac{(\varepsilon w_t)^2}{8}}{w_t\sum_{i=1}^\infty\frac{tp_i}{(tp_i+w_t)^2}+\frac{\varepsilon w_t}{3}}\right).$$
Using Lemma \ref{lem:basic} \ref{item:majo_somme_degueu} we obtain, 
 for $\varepsilon \in (0,1)$, 
\begin{equation}\label{eq:lowerbound2}
    \mathbb{P}(Z^M_{t}\le (1-\eps)w_t)\le  \exp\left(-\frac{\frac{(\varepsilon w_t)^2}{8}}{4w_t+\frac{\varepsilon w_t}{3}}\right)\le \exp\left(-\frac{\varepsilon^2 w_t}{40}\right). 
\end{equation}
Combining \cref{eq:lowerbound1} and \cref{eq:lowerbound2} yields the proposition.
\end{proof}

We deduce from the previous  proposition the following left-tail inequality for $\LIS_n$.

\begin{proposition}\label{prop:concentration2_depoiss}  For all $\eps\in (0,1)$, $n\geq 1$ and every distribution $(p_i)_{i\geq 1}$ we have:

    $$
\mathbb{P}(\LIS_n\le (1-\eps)w_n)\leq 3\exp\left(-\frac{\varepsilon^2(1-\varepsilon)}{40}w_{n}\right).
    $$
\end{proposition}

\begin{proof}
We have 
\begin{eqnarray*}
\mathbb{P}(\LIS_n\le (1-\varepsilon)^2 w_n) & \le & \mathbb{P}(N_{n(1-\varepsilon)}\ge n)+\mathbb{P}(\LIS_n\le (1-\varepsilon)^2 w_n, N_{n(1-\varepsilon)}\le n)\\
& \le & \mathbb{P}(N_{n(1-\varepsilon)}\ge n)+\mathbb{P}(\LIS_{N_{n(1-\varepsilon)}}\le (1-\varepsilon)^2 w_n).
\end{eqnarray*}
Using that $(1-\varepsilon)^2 w_n \le w_{n(1-\varepsilon)}(1-\varepsilon)$ ( Lemma \ref{lem:basic} \ref{ineqw}) and Proposition \ref{prop:concentration2}, we get
$$\mathbb{P}(\LIS_{N_{n(1-\varepsilon)}}\le (1-\varepsilon)^2 w_n)\le \mathbb{P}(\LIS_{N_{n(1-\varepsilon)}}\le (1-\varepsilon) w_{n(1-\varepsilon)})\le 2\exp(-\frac{\varepsilon^2}{40}w_{n(1-\varepsilon)}) \le 2\exp\left(-\frac{\varepsilon^2(1-\varepsilon)}{40}w_{n}\right), $$
using again that $ w_{n(1-\varepsilon)}\ge (1-\varepsilon)w_n $ for the last inequality.
Now, using that $w_n\le n$, we get 
$$\mathbb{P}(N_{n(1-\varepsilon)}\ge n)\le \exp(-\frac{n\varepsilon^2}{40})\le \exp\left(-\frac{w_n\varepsilon^2}{40}\right).$$
\end{proof}

\begin{proof}[Proof of \cref{eq:lower_bound_Theorem}: lower bound in Theorem \ref{th:encadrement}]
We have 
$$\E(\LIS_n)\ge (1-\varepsilon)w_n \mathbb{P}(\LIS_n\ge (1-\varepsilon)w_n).$$
Thus Proposition \ref{prop:concentration2_depoiss} implies that, for all $\varepsilon>0$, 
$$\liminf_{n\to \infty} \frac{\E(\LIS_n)}{w_n}\ge 1-\varepsilon.$$
\end{proof}

\section{Asymptotics of $f_t$}\label{sec:bornes_f}

Theorem \ref{th:encadrement} states that $L_n$ is asymptotically between $\frac{f_n}{2}$ and $f_n$. In this section we try to determine the asymptotic behaviour of the deterministic sequence $f_n$, with strategies that depend on ${\mathbf p}$.

\subsection{Preliminaries: a lower bound from a greedy algorithm}\label{sec:greedy}

Neither Theorem \ref{th:encadrement} nor its proof gives a constructive way to find a long increasing subsequence. This section partially addresses the problem by presenting a simple greedy algorithm which, for a large class of distributions, produces a long increasing subsequence of the right asymptotic order of magnitude.
Besides, this will provide estimates which will be used in the proofs of Proposition \ref{prop:lighttail}  and Proposition \ref{prop:heavytail}.

In order to explain the construction, let us first assume that $(p_i)_{i\ge 1}$ is non-increasing. Fix $n\ge 1$ and  consider the increasing subsequence $(E_1,E_2,\dots, E_{R_n})$ defined as follows:
\begin{itemize}
    \item $E_1$ is the first appearance of $1$ in the sequence $(X_k)_{k\ge 1}$, namely, 
    $
    E_1=\min \set{k\geq 1; X_k=1}
    $.
    \item $E_2$ is the first appearance of $2$ after $E_1$, namely, 
    $
        E_2=\min \set{k\geq E_1; X_k=2}
    $.
    \item ... and so on for any $i\ge 1$.
\end{itemize}
Let $R_n$ be the largest index such that $E_{R_{n}}\le n$.
The sequence $(X_{E_1},X_{E_2},\dots, X_{E_{R_n}})=(1,2,\ldots,R_n)$ is thus an increasing subsequence of $(X_1,\dots , X_n)$. It is clear that $E_i-E_{i-1}$ is a geometric random variable with mean $1/p_i$, so we expect that the length $R_n$ of the above subsequence satisfies $R_n \approx r_n$ where\footnote{Recall that ${\mathbf p}$ has a infinite support, therefore $r_n$ is well defined.} 
$$
r_n=\sup\{r\geq 1, \sum_{i=1}^{r}\frac{1}{p_i}\le n\}.
$$
This construction looks good if $\mathbf{p}$ is non-increasing but may give a very poor result in the case of a general $\mathbf{p}$ (imagine if $p_1$ is very small, or even zero). It is very easy to fix this by reordering the $p_i$'s and this is the purpose of the general construction that we give now:
\begin{enumerate}
    \item For $i\geq 1$, let $p^{\downarrow}_i$ be the $i$-th greatest element of the sequence $\mathbf{p}$ (with multiplicity). In other words $\mathbf{p}^{\downarrow}:=(p^{\downarrow}_i)_{i\ge 1}$ is the non-increasing rearrangement of the sequence $\mathbf{p}$.
    \item For $n\geq 1$ let $r_n$ be uniquely defined by
    \begin{equation}\label{eq:rn}
        r_n=\sup\{r\geq 1, \sum_{i=1}^{r}\frac{1}{p^{\downarrow}_i}\le n\}.
    \end{equation}
    \item Let $i^{(n)}_1<i^{(n)}_2<\dots <i^{(n)}_{r_n}$ be such that
    $$
    \set{p_{i^{(n)}_1},p_{i^{(n)}_2},\dots, p_{i^{(n)}_{r_n}}}=
     \set{p^{\downarrow}_{1},p^{\downarrow}_{2},\dots, p^{\downarrow}_{r_n}}.
    $$
    \item Then apply the construction of the non-increasing case:
 \begin{itemize}
    \item $E_1$ is the first appearance of $i^{(n)}_1$ in the sequence $(X_k)_{k\ge 1}$, namely,  
    $
    E_1=\min \set{k\geq 1; X_k=i^{(n)}_1}
    $.
    \item $E_2$ is the first appearance of $i^{(n)}_2$ after $E_1$, namely, 
    $
        E_2=\min \set{k\geq E_1; X_k=i^{(n)}_2}
    $.
    \item ... and so on for any $i\le r_n$.
\end{itemize}  
\item Finally return $(X_{E_1},\dots,X_{E_{ R_n}})$ where $R_n$ is the largest index such that we have that $E_{R_{n}}\le n$ (note that by construction  $R_n\le r_n$).
\end{enumerate}
This algorithm returns an increasing subsequence of $(X_1,\ldots,X_n)$ of length $R_n$, implying in particular that $L_n\ge R_n$.
Observe that this algorithm is greedy for fixed $n$ but requires to know $n$ in advance in order to determine $i^{(n)}_1,i^{(n)}_2,\dots $. In particular the outputs of the algorithm for successive integers $n,n+1$ may differ significantly and the sequence $(R_n)_{n\ge 1}$ is not in general nondecreasing.
Using this strategy, we obtain the following lower bound for $L_n$.
\begin{proposition}\label{prop:R_n_r_n}
For every distribution $\mathbf{p}$, we have 
$$\liminf_n \frac{\LIS_n}{r_n}\ge 1\qquad \mbox{a.s.}$$
where $(r_n)_{n \ge 1}$ is given by \cref{eq:rn}.
\end{proposition}

\begin{proof} 
Let $\lambda\in (0,1)$ and $n\in \N$. With the notation introduced above, observe that $R_n<\lfloor \lambda r_n\rfloor $ if and only if $E_{\lfloor\lambda r_n\rfloor}>n$.
Since $L_n\ge R_n$, this yields
$$\mathbb{P}(L_n<\lfloor \lambda r_n\rfloor)\le \mathbb{P}(R_n<\lfloor \lambda r_n\rfloor)=\mathbb{P}(E_{\lfloor\lambda r_n\rfloor}>n).$$
Recalling that $E_{\lfloor\lambda r_n\rfloor}$ is the sum of $\lfloor\lambda r_n\rfloor$ independent geometric variables with respective parameters $(p^{\downarrow}_i)_{i\le \lfloor\lambda r_n\rfloor} $, we get
\begin{eqnarray*}
\mathbb{P}(L_n<\lfloor \lambda r_n\rfloor)  
& = &\P\left(E_{\lfloor\lambda r_n\rfloor}-\sum_{i=1}^{\lfloor \lambda r_n\rfloor}\frac{1}{p^{\downarrow}_i} >n-\sum_{i=1}^{\lfloor \lambda r_n\rfloor}\frac{1}{p^{\downarrow}_i}\right)\\
& \leq &\P\left(E_{\lfloor\lambda r_n\rfloor}-\E[E_{\lfloor\lambda r_n\rfloor}] >\sum_{i=1}^{r_n} \frac{1}{p^{\downarrow}_i} -\sum_{i=1}^{\lfloor \lambda r_n\rfloor}\frac{1}{p^{\downarrow}_i}\right)\qquad (\text{using }\sum_{i=1}^{r_n}1/p^{\downarrow}_i\leq n)\\
& \le & \P\left(E_{\lfloor\lambda r_n\rfloor}-\E[E_{\lfloor\lambda r_n\rfloor}]>\sum_{i={\lfloor \lambda r_n\rfloor}+1}^{r_n}\frac{1}{p^{\downarrow}_i}\right)\\
& \le & \frac{\mbox{Var}(E_{\lfloor\lambda r_n\rfloor})}{ \left(\sum_{i={\lfloor \lambda r_n\rfloor}+1}^{r_n}\frac{1}{p^{\downarrow}_i}\right)^2}\\
& \le & \frac{\sum_{i=1}^{\lfloor \lambda r_n\rfloor}\frac{1}{p_i^{\downarrow^2}}}{ \left(r_n(1-\lambda)\frac{1}{p^{\downarrow}_{\lfloor \lambda r_n\rfloor}}\right)^2}
\le  \frac{\lambda}{(1-\lambda)^2 r_n}.
\end{eqnarray*}
 Let us note that $(r_n)_{n \ge 1}$ is a non-decreasing integer sequence which tends to infinity and such that $r_{n+1}-r_{n}\le 1$. Considering $(r_{m_n})_{n\ge 1}$ such that $r_{m_n}=n^2$, the Borel-Cantelli Lemma yields that $L_{m_n}\ge \lfloor\lambda r_{m_n}\rfloor$ a.s. for $n$ large enough.  Using  the same trick as in Section \ref{sec:Cv_ps_Theorem}, we easily deduce 
$$\liminf_m \frac{\LIS_{m}}{r_{m}}\ge 1\qquad \mbox{a.s.}$$
\end{proof}

\subsection{Bounds for $f_t$}\label{sec:bounds_ft}

We have proved that, 
$$\liminf_n \frac{\LIS_n}{r_n}\ge 1 \qquad \mbox{ and } \qquad \limsup_n \frac{\LIS_n}{f_n}\le 1 \qquad \mbox{a.s.}$$
So it is natural to try to find bounds on these two quantities for various  distributions $\mathbf{p}$.

\begin{lemma}\label{lem:Encadrement_f_n}
Assume that  ${\bf p}$ is eventually non-increasing.  Let $p:[0,\infty)\rightarrow (0,\infty)$ be an eventually non-increasing continuous function such that $p(n)=p_n$ for $n\in \N$ and let $\bar{F}(x)=\int_{x}^\infty p(u)du$. Then, for $t$ large enough, there exist  unique $\mu_t>0$ and $\nu_t>0$ such that
\begin{equation}\label{eq:mu_beta}
\frac{\mu_t}{p(\mu_t)}=t \qquad \mbox{ and } \qquad \frac{\nu_t^2}{\bar{F}(\nu_t)}=t.
\end{equation}
Moreover, for $n$ large enough, we have
\begin{equation}\label{eq:Encadrement_f_n}
\lfloor\mu_n\rfloor \le r_n \text{ and } f_n \le 3\nu_n.
\end{equation}
\end{lemma}

\begin{proof}
First note that, for $n$ large enough $\mu_t$ and $\nu_t$ are well defined since the functions $x\mapsto \frac{x}{p(x)}$ and $x\mapsto\frac{x^2}{\bar{F}(x)}$ are both continuous, eventually increasing and tend to infinity. Note also that we necessarily have $\lim_{t\to \infty} \mu_t=\lim_{t\to \infty} \nu_t=\infty$. 
We first prove the first inequality in \cref{eq:Encadrement_f_n}. Since $\mathbf{p}$ is eventually non-increasing and $\mu_n$ tends to infinity, we get, for $n$ large enough,  $p(\mu_n)\le p_i$ for all $i\le \mu_n$.
Thus, we deduce
$$\sum_{i=1}^{\lfloor\mu_n\rfloor}\frac{1}{p_i}\le \frac{\lfloor\mu_n\rfloor}{p(\mu_n)}\le \frac{\mu_n}{p(\mu_n)}= n.$$
So by definition of $r_n$, we get $r_n\ge \lfloor\mu_n\rfloor$.

Let us prove the second inequality in \cref{eq:Encadrement_f_n}.
The function $\frac{\bar{F}(x)}{{x}}$ is continuous, decreasing on $(0,\infty)$ from infinity to   0. Let $q$ be its inverse function, i.e., for all $a>0$,
$\frac{\bar{F}(q(a))}{{q(a)}}=a$ .
Then 
$$\sum_{i\geq 0}\frac{p_i}{p_i+a}\le  \sum_{i\leq q(a)}\frac{p_i}{p_i+a}+ \sum_{i> q(a)}\frac{p_i}{p_i+a}
\leq q(a)+\sum_{i> q(a)}\frac{p_i}{a}
= q(a)+\frac{\bar{F}(q(a))}{a}=2q(a).$$
Thus
$f_n \le  \inf_{a\geq 0}\left\{an+2q(a)\right\}$. Taking $a=\nu_n/n$, we obtain 
\begin{equation*}
    f_n \le  \nu_n+2q(\nu_n/n).
\end{equation*}
By definition of $\nu_n$ we have that 
$\nu_n/n=\bar{F}(\nu_n)/\nu_n$, so we get
$$q(\nu_n/n)=q(\bar{F}(\nu_n)/\nu_n)=\nu_n$$
and so $f_n\le 3 \nu_n$.
\end{proof}

We now show that for light-tailed distributions $r_n$ and $f_n$ are of the same order, implying that our lower and upper bounds for $\LIS_n$ are also of the same order. 
The proposition below applies for example for the geometric and Poisson distributions.

\begin{proposition}[Asymptotics for  sub-exponential  distributions]\label{prop:f_sous_exp} Assume that   
 $\limsup_i p_{i+1}/p_i<1$. Then, when $n\to +\infty$,
 $$
w_n\sim f_n\sim r_n\sim\mu_n.
 $$
 (For convenience we recall that $w_n,f_n,r_n,\mu_n$ are respectively  defined in \cref{eq:wn,eq:fn,eq:mu_beta,eq:rn}.)
\end{proposition}
Note that Proposition \ref{prop:f_sous_exp} combined with Theorem \ref{th:encadrement} implies that $\frac{L_n}{f_n}\to 1$ for all sub-exponential distributions.
\begin{proof}
        First, combining  Theorem \ref{th:encadrement} with  Proposition \ref{prop:R_n_r_n}, we get $\liminf_n \frac{f_n}{r_n}\geq 1$. Moreover, with Lemma \ref{lem:Encadrement_f_n}, we know that 
        $\liminf_n \frac{f_n}{\mu_n}\ge 1$ (the assumption  $\limsup_i p_{i+1}/p_i<1$ ensures that $(p_i)_{i\ge 1}$ is eventually non-increasing).
        
        Let us show that $\limsup_n \frac{f_n}{\mu_n}\le 1$. Let $\eps>0$ and $n\in\N$. We write (recall that $g_n$ was defined in \cref{eq_def_g})
        \begin{align*}
            g_n(\eps\frac{\mu_n}{n})&=\eps\mu_n+\sum_{i=1}^{\lfloor\mu_n\rfloor}\frac{p_i}{p_i+\eps\frac{\mu_n}{n}}+\sum_{i=\lfloor \mu_n\rfloor+1}^{+\infty}\frac{p_i}{p_i+\eps\frac{\mu_n}{n}}\\
                &\le \eps\mu_n +\sum_{i=1}^{\lfloor\mu_n\rfloor}1+\sum_{i=\lfloor \mu_n\rfloor+1}^{+\infty}\frac{n}{\eps\mu_n}p_i\\
                &\le\eps\mu_n+\mu_n+\frac{n}{\eps\mu_n}\sum_{i>\mu_n}p_i.
        \end{align*}
        Now, since $\limsup_i \frac{p_{i+1}}{p_i}<1$, there exists $0<\lambda<1$ and $j\in\N$ such that for all $i\geq j$ and $k\geq0$, $p_{i+k}\le p_i\lambda^k$. Thus for all $i\in\N$
        $$
        \sum_{k>i}p_k\le p_i\sum_{k=0}^{+\infty}\lambda^k=p_i\frac{1}{1-\lambda}.
        $$
        It yields that 
        $$
f_n\le g_n(\eps\frac{\mu_n}{n})\le (1+\eps)\mu_n+ \frac{n}{\eps\mu_n}p(\mu_n)\frac{1}{1-\lambda}.
        $$
        Since $\frac{np(\mu_n)}{\mu_n}=1$ we finally have
        $$
f_n\le (1+\eps)\mu_n+\frac{1}{(1-\lambda)\eps},
        $$
        which implies
        $
\limsup_n \frac{f_n}{\mu_n}\le 1
$
(recall that $\mu_n\to +\infty$).         
        Now, according to  Lemma \ref{lem:Encadrement_f_n}, for $n$ large enough $\lfloor\mu_n\rfloor\le r_n$ so 
        $$
\limsup_n \frac{f_n}{r_n}\le \limsup_n \frac{f_n}{\mu_n}\le 1,
        $$
        thus
        $$
\lim_{n\to+\infty}\frac{f_n}{r_n}=\lim_{n\to+\infty}\frac{f_n}{\mu_n}=1
        $$
      \emph{i.e.}  $f_n\sim\mu_n\sim r_n$. Since, $w_n\le f_n$, it remains to find a lower bound of $w_n$  of the same order. Let $\varepsilon>0$ and let us prove that, for $n$ large enough, $w_n\ge \mu_n(1-\varepsilon):=t_n$.
        For $n\ge 0$, put $$h_n(x):=\sum_{i\ge 1} \frac{np_i}{np_i+x}$$
        so that $h_n(w_n)=w_n$. We show that $h_n(t_n)\ge t_n$.
        Since $h_n$ is a non-increasing function, this would imply that $w_n\ge t_n$ and conclude the proof.
        Let $\eta>0$ such that $(1-\eta)/(1+\eta)=(1-\varepsilon)$ and let 
        $$I_n:=\{i\ge 1, \frac{t_n}{np_i}\le \eta\}.$$
We have 
$$h_n(t_n)\ge \sum_{i\in I_n} \frac{1}{1+\frac{t_n}{np_i}}\ge |I_n|\frac{1}{1+\eta}.$$
We need now to find a lower bound on the cardinal of $I_n$. We claim, that for $n$ large enough, 
$$|I_n|\ge (1-\eta)\mu_n.$$
Indeed recalling that $p(\mu_n)=\mu_n/n$ We have
$$i\in I_n \Longleftrightarrow p_i\ge \frac{t_n}{\eta n} \Longleftrightarrow p(\mu_n)\le \frac{\eta}{(1-\varepsilon)}p_i.$$
Using that $\mathbf{p}$ is sub-exponential, we have, for $n$ large enough and $i\le (1-\eta)\mu_n $, 
$$p(\mu_n)\le \lambda^{\eta \mu_n} p((1-\eta)\mu_n)\le \lambda^{\eta \mu_n} p_i$$
and so $i\in I_n$ if $n$ is large enough such that $\lambda^{\eta \mu_n}\le \eta/(1-\varepsilon)$.
Hence, we get
$$h_n(t_n)\ge \frac{1-\eta}{1+\eta}\mu_n=t_n.$$
    \end{proof}

We now give a lemma comparing $(\mu_t)_{t>0}$ and $(\nu_t)_{t>0}$ for heavy tailed distribution regularly varying with index strictly larger than $1$. This lemma will be useful for the proof of  Proposition \ref{prop:heavytail} below.
\begin{lemma}\label{lem:borne_nu_mu}
If $p_i= \ell(i)i^{-\gamma}$ with $\gamma>1$ and $\ell$ slowly varying and $\mathbf{p}$ is eventually non-increasing, then, for all $c>1$, for $t$ large enough $ \nu_t \le \mu_{ct/(\gamma-1)}$.
\end{lemma}

\begin{proof}
If $p_i= \ell(i)i^{-\gamma}$, then $\bar{F}(i):= \sum_{j>i}p_j\sim \frac{1}{\gamma-1}\ell(i)i^{-\gamma+1}$ (see \cite{feller1991},  Section VIII.9, Theorem 1). In particular, for any $c>1$ and $x$ large enough, 
$$\bar{F}(x)\le  \frac{c}{\gamma-1}x p(x).$$
Thus, for $t$ large enough
$$\nu_t^2=t\bar{F}(\nu_t)\le t\frac{c }{\gamma-1}\nu_t p(\nu_t)$$
and so 
$$\frac{ \nu_t}{p( \nu_t )}\le \frac{tc }{\gamma-1}= \frac{ \mu_{tc/(\gamma-1)} }{p( \mu_{tc/(\gamma-1)})}.$$
We use again that the function $x\rightarrow x/p(x)$ is increasing to conclude.
\end{proof}

\subsection{Some explicit estimates}\label{sec:explicit}

Recall that $\mu_n$ was defined in \cref{eq:mu_beta} as the solution of 
$\frac{\mu_n}{p(\mu_n)}=n$.  The combination of 
Lemma \ref{lem:Encadrement_f_n}, Proposition \ref{prop:f_sous_exp} and Lemma \ref{lem:borne_nu_mu} allow to translate bounds for $\mu_n$ into  bounds for $f_n$, for various distributions.
The purpose of the next lemma is therefore to estimate $\mu_n$ in several cases.

\begin{lemma}\label{lem:calcul} 
 \begin{itemize}
    \item If ${\bf p}$ is a geometric distribution $p_i:=p (1-p)^{i-1}$, then  
    \begin{equation}\label{eq:mu_geometrique}
    \mu_n\sim \frac{\log n}{|\log (1-p)|}.
    \end{equation}
    \item If  ${\bf p}$ is a Poisson distribution $p_i:=\frac{e^{-\lambda}\lambda^i }{i!}$, then 
    \begin{equation}\label{eq:mu_poisson}
    \mu_n\sim \frac{\log n}{|\log \log n|}.
    \end{equation}
      \item If  ${\bf p}$ is eventually non-increasing and $p_i\sim c i^{-\beta}(\log i)^\gamma$, $\beta>1$, $\gamma \in \R$, then 
      $$
      \mu_n\sim \left( c n\left(\frac{\log n}{(1+\beta)}\right)^\gamma \right)^{1/(1+\beta)}.
      $$
 \end{itemize}
 \end{lemma}
 \begin{proof}
Recall that  $p:[0,\infty)\rightarrow (0,\infty)$ is an eventually non-increasing continuous function such that $p(i)=p_i$ for $n\in \N$. It follows that $\lim_{n\to +\infty}\mu_n =+\infty$.
\begin{itemize}
\item  Assume that $p_i=p(1-p)^{i-1}$ and define the function $p$ for  $x\ge 1$ by $p(x)=p(1-p)^{x-1}$. Then
 $$
 \frac{\mu_n}{n}=p(\mu_n)=p(1-p)^{\mu_n-1}.
 $$
 Thus we have
 $$
 \log \mu_n - \mu_n\log(1-p) =\log n+\log\left(\frac{p}{1-p}\right).
 $$
 Now, since $\log{\mu_n}=o(\mu_n)$, we get $\mu_n|\log(1-p)|\sim \log n$, hence \cref{eq:mu_geometrique}.
 \item Consider now the distribution $p_i=\frac{e^{-\lambda}\lambda^i}{i!}$ and  for  $x\ge 1$ set $p(x)= e^{-\lambda}\frac{\lambda^x}{\Gamma(x+1)}$. The function $p$ is increasing on $[1,x_\lambda]$ and decreasing on $[x_\lambda,+\infty)$ for some $x_\lambda >1$.
Since $\mu_n\to+\infty$, there exists $n(\lambda)\in\N$ such that for all $n\geq n(\lambda)$, $\mu_n\geq x_\lambda$. 
For such $n$,
 $$
 p(\mu_n)=\frac{\mu_n}{n}=\frac{e^{-\lambda}\lambda^{\mu_n}}{\Gamma(\mu_n+1)}.
 $$
 Now, since $\Gamma(\mu_n+1)=\mu_n\Gamma(\mu_n)$ we get
 $$
 \mu_n^2\frac{\Gamma(\mu_n)}{\lambda^{\mu_n}}=\frac{n}{e^\lambda}.
 $$
 Thanks to the Stirling formula we can write:
 $$
 \mu_n^2\times \left(\frac{\mu_n}{\lambda e}
 \right)^{\mu_n}\sqrt{2\pi \mu_n}\sim \frac{n}{e^\lambda}.
 $$
 Now by taking the logarithm and using that $\log \mu_n$ and $ \mu_n$ are negligible compared to $\mu_n \log \mu_n$, we get $ \mu_n\log(\mu_n)\sim \log(n)$. Applying the logarithm a last time, we  then obtain  $\log(\mu_n)\sim \log(\log(n))$, hence \cref{eq:mu_poisson}.
 \item Let $\beta>1$, $\gamma\in \R$ and assume that $p_i\sim c(\log i)^\gamma i^{-\beta}$,  so $p(\mu_n)\sim c(\log \mu_n)^\gamma\mu_n^{-\beta}$ for $n$ large enough. We get
 $$n=\frac{\mu_n}{p(\mu_n)}\sim \frac{\mu_n^{1+\beta}}{c(\log \mu_n)^\gamma}.$$ 
 Taking the logarithm, we get 
 $$\log n \sim (1+\beta) \log \mu_n.$$
 Plugging this in the previous equation, we get 
 $$\mu_n^{1+\beta}\sim c n\left(\frac{\log n}{(1+\beta)}\right)^\gamma.$$ 
 which yields the result.
\end{itemize}
 \end{proof}
 
 \begin{proof}[Proof of Proposition \ref{prop:lighttail}] Using  Theorem \ref{th:encadrement}, we know that 
 $$\liminf_n \frac{L_n}{w_n}\ge 1 \qquad \limsup_n \frac{L_n}{f_n}\le 1 \qquad a.s.$$
 Moreover, if $\mathbf{p}$ is a geometric or a Poisson distribution, thanks to 
 Proposition \ref{prop:f_sous_exp}, we know that 
 $f_n\sim w_n \sim \mu_n$. 
 Then we estimate $\mu_n$ using  Lemma \ref{lem:calcul} in order to deduce Proposition \ref{prop:lighttail}.
 
 \end{proof}

\begin{proof}[Proof of Proposition \ref{prop:heavytail}]
Recall that $\mu_n$ and $\nu_n$ are defined by \cref{eq:mu_beta}.
Combining Theorem \ref{th:encadrement} and Proposition \ref{prop:R_n_r_n} we obtain
\begin{equation}\label{eq:encadrement_mu_nu_inf}
    \liminf_n \frac{L_n}{\mu_n}=  \liminf_n \frac{L_n}{r_n}\times \frac{r_n}{\mu_n} \geq 1 \times 1 >0. 
\end{equation}
Besides, combining Theorem \ref{th:encadrement}, and Lemma \ref{lem:Encadrement_f_n} yields  
\begin{equation}\label{eq:encadrement_mu_nu_sup}
\limsup_n \frac{L_n}{\nu_n}= 
\limsup_n \frac{L_n}{f_n} \times  \frac{f_n}{\nu_n}
\leq 1 \times 3 <+\infty.
\end{equation}

In the case $p_i\sim c (\log i)^\gamma i^{-\beta}$ with $\beta >1$, we deduce \ref{PropHT_i} in  Proposition \ref{prop:heavytail}  using the asymptotic for $\mu_n$ computed in  Lemma \ref{lem:calcul}  and the upper bound for $\nu_n$ given in Lemma \ref{lem:borne_nu_mu}.

Let us now consider the case \ref{PropHT_ii} of Proposition \ref{prop:heavytail}, namely $p_i\sim c (\log i)^\gamma i^{-1}$ with $\gamma <-1$. 
Then,  we have, for some $c_1>0$,
$$ \bar{F}(x)= \sum_{i> x} p_i \sim c_1(\log x)^{\gamma+1}.$$
As $\frac{\nu_n^2}{\bar{F}(\nu_n)}=n$ this implies that  $$\nu_n\sim c_2\sqrt{n(\log n)^{\gamma+1}}$$
for some $c_2 >0$. 
This asymptotic combined with \cref{eq:encadrement_mu_nu_sup} yields the upper bound 
$$\limsup_n \frac{L_n}{\sqrt{n(\log n)^{\gamma+1}}}<\infty.$$
An easy computation shows that $\mu_n\sim c_3{\sqrt{n(\log n)^{\gamma}}}$, so we will not obtain  the desired lower bound for $L_n$ just using that asymptotically $L_n\ge \mu_n$. We should make here more precise computation to find a lower bound for $f_n$ (and so for $w_n$).

Recall that $p(x)$ is eventually non-increasing and such that $p(i)=i$ for every  integer $i$. For $a$ close to $0$, let $u(a)$ such that $a=p(u(a))$. Thus
$$
a \stackrel{a\to 0}{\sim} c (\log u(a))^\gamma u(a)^{-1}
$$
which implies 
$$ 
u(a)\stackrel{a\to 0}{\sim} c |\log a|^\gamma a^{-1}.
$$
Using that $\mathbf{p}$ is eventually decreasing, we have, for $a$ close to $0$, $p_i\le a$ if $i\ge u(a)$. So we get 
$$\sum_{i=1}^\infty \frac{p_i}{p_i+a} \ge \sum_{i\ge u(a)} \frac{p_i}{a+a} \sim \frac{\bar{F}(u(a))}{2a} \sim  \frac{c_1|\log a|^{\gamma+1}}{2a}.$$
Hence, we get, for $n$ large enough 
\begin{equation*}
    f_n= \inf_{a\geq 0}\left\{ an+\sum_{i\geq 0}\frac{p_i}{p_i+a}\right\}\ge \inf_{a\geq 0}\left\{ an+\frac{c_1|\log a|^{\gamma+1}}{4a }\right\}\sim c_3\sqrt{n(\log n)^{\gamma+1}}.
\end{equation*}
Using Theorem \ref{th:encadrement} and that $w_n\ge f_n/2$, this yields the expected lower bound:
$$\liminf_n \frac{L_n}{\sqrt{n(\log n)^{\gamma+1}}}>0.$$

\end{proof}

\begin{remark}\label{rem:indice1}
In the case $p_i\sim c (\log i)^\gamma i^{-1}$ with $\gamma <-1$, we can also compute the order of $r_n$. Indeed, up to constant,
$\sum \frac{1}{p_i}\sim c\sum i (\log i)^{-\gamma} \sim c_4k^2 (\log k)^{-\gamma}$
And so $r_n$ is such that $c_4r_n^2 (\log r_n)^{-\gamma} \sim n$ which yields $r_n\sim c_5\sqrt{n (\log n)^\gamma}$. Hence, we see that $L_n\gg r_n$ which shows that in this case, the greedy algorithm described in Section \ref{sec:greedy} does not give a subsequence of the right order.
\end{remark}

\section{Case of a generic distribution: proof of Theorem  \ref{th:distribution_generique}}\label{sec:distribution_generique}

We explain in this section how Theorem \ref{th:encadrement}, dealing with  discrete distributions on integers, can in fact be extended to any discrete distribution, even if its support cannot be enumerated in increasing order, as in the cases of $\mathbb{Z}$ or $\mathbb{Q}$.
Let us first observe that the proof of $\LIS_n/\mathbb{E}(\LIS_n)\to 1$ (\emph{i.e} \cref{eq:Cv_ps_Theorem}) is valid for any distribution $\rho$. Therefore, in order to extend the Theorem \ref{th:encadrement} to arbitrary distributions  we just need to control $\E(L_n)$.

\begin{proof}[Proof of Theorem \ref{th:distribution_generique}]
We begin by proving item (i) in Theorem \ref{th:distribution_generique}.
Let $\rho$ be a discrete distribution on a countable set $S:=\{i_1,i_2,\ldots\}$. Let $(X_i)_{i\ge 1}$ be an i.i.d. sequence with law $\rho$ and for $M,n\ge 1$, set
$$D_n^M:=\{1\le i \le n, X_i\in \{i_1,\ldots,i_M\}\}.$$
We define $L_n^M:=L((X_i)_{i\in D_n^M})$, the length of the longest increasing subsequence with values in $\{i_1,\dots,i_M\}$.
Since the set $\{i_1,\dots,i_M\}$ is finite, it can be enumerated in the increasing order and so we can apply Theorem \ref{th:poisson_domination} : if $N_t$ is a Poisson random variable with parameter $t$ independent of the sequence $(X_i)_{i\ge 1}$, then, for any $\alpha>0$, we have the stochastic domination :
\begin{equation*}
    L^M_{N_t} \leq \mathcal{H}^{(\alpha)}_t+ \sum_{j=1}^M I^{(\alpha)}_j,
    \end{equation*}
where $\mathcal{H}^{(\alpha)}_t$ is a Poisson random variable with parameter $\alpha t$ and  $(I^{(\alpha)}_j)_{j \geq 1}$ is a sequence of independent random variables with $\mathrm{Bernoulli}\left(\frac{\rho(\{i_j\})}{\alpha+\rho(\{i_j\})}\right)$ distribution. Thus, for any $\alpha>0$
\begin{equation*}
   \E( L^M_{N_t} )\leq \alpha t + \sum_{j=1}^M \frac{\rho(\{i_j\})}{\alpha+\rho(\{i_j\})} \le \alpha t + \sum_{i\in S} \frac{\rho(\{i\})}{\alpha+\rho(\{i\})}.
    \end{equation*}
Optimizing over $\alpha$, we get     $\E( L^M_{N_t} )\le f_t$ and so, by monotone convergence, $\E( L_{N_t} )\le f_t$. We conclude as in Section \ref{SectionUpper} that 
$$\limsup_n \frac{\E(L_n)}{f_n}\le 1.$$

 \newcommand{\alfa}{\widetilde{\alpha}_t}
Regarding the lower bound, using again Theorem \ref{th:poisson_domination}, we get the stochastic lower bound
\begin{equation*}
  L_{N_t} \geq   L^M_{N_t} \geq \min(\mathcal{H}^{(\alpha)}_t, \sum_{j=1}^M J^{(\alpha)}_j).
    \end{equation*}
Let $w_t$ the unique positive solution of the equation  
\begin{equation*}
w_t=\sum_{i\in S} \frac{t\rho(\{i\})}{t\rho(\{i\})+w_t}
\end{equation*}
and take   $\alfa=w_t/t$. We get
 \begin{equation*}
    \mathbb{P}(L_{N_t}<(1-\eps)w_t)\leq \mathbb{P}(\mathcal{H}^{(\alfa)}_t\le (1-\eps)w_t)+ \mathbb{P}(Z^M_{t}\le (1-\eps)w_t),
\end{equation*}
where  $Z^M_{t}:=\sum_{i=1}^M J^{(\alfa)}_i$. One can check that \cref{eq:lowerbound1} and \cref{eq:lowerbound2} still hold in this case if $M$ is chosen large enough and thus imply that 
  $$
\mathbb{P}(L_{N_t}<(1-\eps)w_t)\leq 2\exp\left(-\frac{\varepsilon^2 w_t}{40}\right),
    $$
as obtained in Proposition \ref{prop:concentration2} in the case of a distribution on the integers. We conclude, as in Section \ref{SectionLower}, that 
$$\liminf_n \frac{\E(L_n)}{w_n}\ge 1.$$

\medskip

Let us now prove item (ii) in Theorem \ref{th:distribution_generique}. Let $\rho$ be a probability measure on $\R$ and write $\rho=\rho_1+\rho_2$ where $\rho_1$ is atomless and $\rho_2$ is purely discrete. Let $S$ be the (countable) support of $\rho_2$. For $n\ge 1$, define
$$D^1_n:=\{i\le n, X_i\notin S\} \qquad \mbox{ and }\qquad D^2_n:=\{i\le n, X_i\in S\}$$
and set $L_n^1:=L((X_i)_{i\in D^1_n})$ and  $L_n^2:=L((X_i)_{i\in D^2_n})$. Note that, for $k\in \{1,2\}$,  $|D^k_n|$ is a binomial random variable with parameters $(n,\rho_k(\R))$.   Moreover,  since $\rho_1$ is atomless, conditionally on $|D^1_n|$, $L_n^1$ has the law of $\mathcal{L}_{|D^1_n|}$, the length of the longest increasing subsequence in a uniform permutation of size  $|D^1_n|$. 
Besides, we have
$$L_n^1\le L_n\le L_n^1+L_n^2.$$
On the one hand, using (i) of Theorem \ref{th:distribution_generique} and the fact that $f_n=o(\sqrt{n})$ for any discrete distribution, we get 
$$\lim_{n\to \infty} \frac{L_n^2}{\sqrt{n}}=0 \mbox{ a.s.}$$
On the other hand, using the Ver\v{s}ik-Kerov's theorem on uniform permutations \cite{vershik1977asymptotics}, we obtain
 $$\lim_{n\to \infty} \frac{L_n^1}{2\sqrt{|D_n^1|}}=1 \mbox{ a.s.}$$
as soon as $|D_n^1|$ converges to infinity almost surely. We conclude using the fact that, if $\rho_1(\R)>0$, then $|D_n^1|\sim n\rho_1(\R)$ almost surely.

 \end{proof}

\bibliographystyle{alpha}
\bibliography{BiblioLIS}

\end{document}